\input amstex
\documentstyle{amsppt}

\input label.def
\input degt.def

\input epsf
\def\picture#1{\epsffile{#1-bb.eps}}
\def\cpic#1{$\vcenter{\hbox{\picture{#1}}}$}

{\catcode`\@=11
\gdef\proclaimfont@{\sl}
\gdef\subsubheadfont@{\subheadfont@}
}

\loadbold

\def\todo:{{\bf to do:}}

\let\=\B
\def\latin{\emph}
\def\ie{\latin{i.e\.}}
\def\eg{\latin{e.g\.}}
\def\cf.{\latin{cf\.}}
\def\all.{\latin{et al\.}}
\def\etc.{\latin{etc\.}}
\def\via{\latin{via}}

\newhead\subsubsection\subsubsection\endsubhead
\newtheorem{Remark}\Remark\thm\endAmSdef
\speculation\thm\endAmSdef

\def\dash{\item"\hfill--\hfill"}
\def\nodash{\item""}
\def\Dashes{\widestnumber\item{--}\roster}
\def\endDashes{\endroster}

\def\GAP{{\tt GAP}}

\def\mathrm{\roman}

\def\PP{\Bbb P}
\def\Cp#1{\PP^{#1}}
\def\Rp#1{\Cp{#1}_\R}

\def\ls|#1|{\mathopen|#1\mathclose|}

\def\PSL{\operatorname{\text{\sl PSL}}}
\def\SL{\operatorname{\text{\sl SL}}}
\def\MW{\operatorname{\text{\sl MW}}}
\def\NS{\operatorname{\text{\sl NS}}}

\def\depth{\operatorname{dp}}
\def\Sp{\operatorname{Sp}}
\def\Sk{\operatorname{Sk}}
\def\Conv{\operatorname{conv}}
\def\Aut{\operatorname{Aut}}
\def\divides|{\mathbin|}

\let\Ga\alpha
\let\Gb\beta
\let\Gg\gamma
\let\Gs\sigma
\let\Gr\rho

\def\tGa{\tilde\Ga}

\def\tGg{\tilde\Gg}

\let\onto\twoheadrightarrow
\def\scirc{^\circ}
\def\bcirc{^{\circ\circ}}
\let\<\langle
\let\>\rangle
\def\1{^{-1}}

\def\tj{\tilde\jmath}
\def\tx{\tilde x}
\def\ty{\tilde y}

\def\CH{\Cal H}
\def\CK{\Cal K}

\def\CO{\Cal O}
\def\CQ{\Cal Q}
\def\CZ{\Cal Z}
\def\chS{S\spcheck}
\def\X{\Bbb X}
\def\Y{\Bbb Y}

\def\bA{\bold A}
\def\bD{\bold D}
\def\bE{\bold E}
\def\bJ{\bold J}

\def\tA{\tilde\bA}
\def\tD{\tilde\bD}
\def\tE{\tilde\bE}
\def\tJ{\tilde\bJ}
\def\I{\mathrm{I}}
\def\II{\mathrm{II}}
\def\III{\mathrm{III}}
\def\IV{\mathrm{IV}}

\def\tG{\tilde G}
\def\tX{\tilde X}
\def\tF{\tilde F}

\def\fH{\frak H}

\def\IK{\Cal I}

\let\BB=B  
\let\CC=C  
\let\EE=E  
\let\DD=D  
\let\disk\Delta
\let\slope\varkappa
\def\tslope{\tilde\slope}
\def\BM{\frak{Im}} 
\let\AM=A          
\def\AM{\mathrm{A}} 

\let\MG\Gamma            
\def\tMG{\tilde\MG}      
\def\MGi#1#2{\MG_{\!#1}(#2)}
\def\MGd{\MGi1}          
\def\CG#1{\Z_{#1}}       
\def\GDG{\Bbb D}         
\def\DG#1{\GDG_{#1}}     
\def\BG#1{\Bbb B_{#1}}   
\def\SG#1{\Bbb S_{#1}}   
\def\AG#1{\Bbb A_{#1}}   
\def\FF#1{\Bbb F_{#1}}   
\def\AA{\Bbb G}          
\def\AA{\frak F}         

\def\bm{\frak m}
\let\kk\kappa

\def\stext#1{^{\botsmash{\operatorname{#1}}}}
\def\sproj{\stext{proj}}
\def\saff{\stext{afn}}
\def\piproj#1{\pi\sproj_{#1}}
\def\piaff#1{\pi\saff_{#1}}
\def\pibul#1{\pi^\bullet_{#1}}

\def\inserthyphen{\ifcat\next a-\fi\ignorespaces}
\let\BLACK\bullet
\let\WHITE\circ
\def\CROSS{\vcenter{\hbox{$\scriptstyle\mathord\times$}}}
\let\STAR*
\def\TRIANG{\vcenter{\hbox{$\scriptstyle\mathord\vartriangle$}}}
\def\pblack-{$\BLACK$\futurelet\next\inserthyphen}
\def\pwhite-{$\WHITE$\futurelet\next\inserthyphen}
\def\pcross-{$\CROSS$\futurelet\next\inserthyphen}
\def\pstar-{$\STAR$\futurelet\next\inserthyphen}
\def\ptriang-{$\TRIANG$\futurelet\next\inserthyphen}
\def\black{\protect\pblack}
\def\white{\protect\pwhite}

\topmatter

\title
Dihedral coverings of trigonal curves
\endtitle

\author
Alex Degtyarev
\endauthor

\address
Department of Mathematics,\endgraf\nobreak
Bilkent University,\endgraf\nobreak
06800 Ankara, Turkey
\endaddress

\email
degt\@fen.bilkent.edu.tr
\endemail

\abstract
We classify and study
trigonal curves in Hirzebruch surfaces
admitting dihedral Galois coverings. As a consequence, we obtain
certain restrictions on the fundamental group of a plane
curve~$D$ with a singular point of multiplicity $(\deg D-3)$.
\endabstract



\endtopmatter

\document

\section{Introduction\label{S.intro}}

\subsection{Motivation}
This paper begins a systematic study of the fundamental group of
a trigonal curve in a geometrically ruled rational surface. The
main tools used are the braid monodromy, Zariski--van Kampen
theorem, and arithmetic properties of the modular group
$\MG:=\PSL(2,\Z)$.

Originally, my interest in trigonal curves was motivated by my
attempts to compute the fundamental groups of plane
curves, a problem that was first posed by
O.~Zariski~\cite{Zariski.group}, \cite{Zariski} in 1930
and that has since
been a subject of intensive research by a number of
mathematicians. Given such a curve $\DD\subset\Cp2$, one can blow
up a singular point~$P$ and, by a sequence of elementary
transformations, convert~$\DD$ to a curve $\CC\subset\Sigma_d$ in
a Hirzebruch surface. If the point~$P$ is of multiplicity
$(\deg\DD-3)$, then~$\CC$ is a trigonal curve. The fundamental
group of~$\DD$ is closely related, although not
isomorphic, to that of~$\CC$
(see Subsections~\ref{s.generalized}
and~\ref{proof.plane} for details),
and any information on the former can shed a light on
the structure of the latter.
Our principal result in this direction is Theorem~\ref{th.plane}
below. At the same time, it turns out that trigonal curves are of
a certain interest on their own right (for example, as the
ramification loci of elliptic surfaces), and in the framework of
trigonal curves many
statements relating the fundamental group and other geometric
properties take more precise and complete form, see, \eg,
Theorem~\ref{th.Oka}.

\subsection{Principal results}\label{s.results}
Given a trigonal curve $\CC\subset\Sigma_d$ (see
Subsection~\ref{s.curves} for the precise definitions), we will
consider both the \emph{projective} and
\emph{affine fundamental groups}
$$
\piproj\CC:=\pi_1(\Sigma_d\sminus(\CC\cup\EE)),\qquad
\piaff\CC:=\pi_1(\Sigma_d\sminus(\CC\cup\EE\cup F));
$$
here
$\EE\subset\Sigma_d$ is the exceptional section and
$F\subset\Sigma_d$ is a fiber of the ruling
that is nonsingular for~$\CC$. The latter is a cyclic central
extension of the former, and the commutants of the two groups are
equal, see Corollary~\ref{commutants}.
With the usual
abuse of the language, we refer to~$\piproj\CC$ and $\piaff\CC$ as
the (fundamental)
groups of~$\CC$ rather than mentioning explicitly the complement
$\Sigma_d\sminus\ldots$.
Similarly, speaking about a \emph{covering}
of~$\CC$, we mean a covering of~$\Sigma_d$ ramified at
$\CC\cup\EE$.

Our main goal is substantiating the following speculation
(in which Item~\itemref{spec}{spec.universal} is, in fact, a
statement, see Proposition~\ref{G.universal}).

\speculation\label{spec}
\roster
\item\local{spec.bounds}
There do exist certain strict bounds
on the complexity of the fundamental
group of a trigonal curve; they are
due to the fact that the monodromy
group of such a curve is of genus zero.
\item\local{spec.universal}
Any trigonal curve~$\CC$ whose group admits a prescribed quotient
$\piaff\CC\onto G$ is essentially induced from a certain
\emph{universal curve} with this property.
\item\local{spec.geometry}
As a consequence, the existence of a quotient $\piaff\CC\onto G$
as above may imply
certain additional geometric properties of~$\CC$.
\item\local{spec.larger}
In particular, the existence of a quotient $\piaff\CC\onto G$ may
imply the existence of a larger quotient
$\piaff\CC\onto\smash{\tilde G}\onto G$.
\endroster
\par\removelastskip
\endspeculation

As a first step supporting Speculation~\ref{spec}, we discuss the
generalized dihedral quotients of~$\pibul\CC$. Given an abelian
group~$\CQ$, the \emph{\rom(generalized\rom) dihedral group}
$\GDG(\CQ)$ is
the semi-direct product $\CQ\rtimes\CG2$, with
the factor~$\CG2$ acting
on the kernel~$\CQ$ \via~$-\id$.
We use the standard abbreviation $\DG{2n}=\GDG(\CG{n})$ for the
classical dihedral groups; note that, in our notation, the index
refers to the order of the group.

We are interested in the so called \emph{uniform}
dihedral quotients $\pibul\CC\onto\GDG(\CQ)$, see
Definition~\ref{def.uniform}. Roughly, the corresponding covering
is required to have the same ramification behavior over each
irreducible component of~$\CC$. If $\CC$ is irreducible, then each
dihedral quotient of $\pibul\CC$ (respectively, each dihedral
covering of~$\CC$) is uniform, see Proposition~\ref{all.uniform}.

\theorem\label{th.D}
If
the group $\piaff\CC$ of
a nontrivial \rom(see Definition~\ref{def.trivial}\rom)
trigonal curve $\CC\subset\Sigma_d$ admits a
uniform
quotient $\GDG(\CQ)$,
then $\CQ$ is a quotient of one of the
following groups\rom:
$$
\CG2\oplus\CG8,\quad
\CG4\oplus\CG4,\quad
\CG2\oplus\CG6,\quad
\CG3\oplus\CG6,\quad
\CG9,\quad
\CG5\oplus\CG5,\quad
\CG{10},\quad
\CG7.
$$
All quotients of the groups above do appear.
\endtheorem

\theorem\label{th.reducible}
A trigonal curve $\CC\subset\Sigma_d$ is reducible
\rom(respectively,
splits into three
components\rom) if and only
if the group $\piaff\CC$ admits a quotient to
$\DG4=\CG2\oplus\CG2$ \rom(respectively,
to $\GDG(\CG2\oplus\CG2)=\CG2\oplus\CG2\oplus\CG2$\rom).
\endtheorem

\corollary\label{cor.irreducible}
If an
irreducible trigonal curve
admits a $\GDG(\CQ)$-covering, then $\CQ$ is a quotient of
$\CG3\oplus\CG3$, $\CG9$, $\CG5\oplus\CG5$, or $\CG7$.
\endcorollary

Theorem~\ref{th.D} and Corollary~\ref{cor.irreducible} are proved
in Subsection~\ref{proof.D};
Theorem~\ref{th.reducible} is proved in
Subsection~\ref{proof.reducible}.
Note that Theorem~\ref{th.D}, listing a finite set of options, is
in a sharp contrast with the case of hyperelliptic curves in
Hirzebruch surfaces, where each dihedral group~$\DG{2n}$ can
appear as a uniform dihedral quotient of the fundamental group,
see Remark~\ref{rem.BG2} below.

One can notice a certain similarity between trigonal curves and
plane sextics, where the dihedral quotients of the fundamental
groups are also known, see~\cite{degt.Oka}. Although the
particular lists of groups differ, the prime factors appearing in
their orders are the same: one has $p=3$, $5$, $7$, and, for
reducible curves, $p=2$. As another similarity and, at the same
time, an illustration of
Statements~\iref{spec}{spec.geometry},~\ditto{spec.larger},
one has the following almost literal
translation of the stronger version of
Oka's conjecture~\cite{Oka.conjecture} on the Alexander polynomial
of an irreducible plane sextic, see~\cite{degt.Oka}.

\theorem\label{th.Oka}
For an irreducible trigonal curve $\CC\subset\Sigma_d$, the
following four statements are equivalent\rom:
\roster
\item\local{Oka.D}
$\piproj\CC$ factors to
the dihedral group $\DG6\cong\SG3$\rom;
\item\local{Oka.B}
$\piproj\CC$ factors to
the
modular group $\MG\cong\CG2*\CG3$\rom;
\item\local{Oka.Delta}
$t^2-t+1$ divides
the Alexander polynomial $\Delta_\CC(t)$,
see Subsection~\ref{s.Alexander}\rom;
\item\local{Oka.torus}
$\CC$ is of torus type, see Subsection~\ref{s.torus}.
\endroster
\par\removelastskip
\endtheorem

This theorem is proved in Subsection~\ref{proof.Oka}, and its
extension to reducible curves is discussed in
Remark~\ref{rem.Oka}.
A number of other examples
illustrating~\iref{spec}{spec.larger}
are found in
Section~\ref{S.groups}, see, \eg,
Corollaries~\ref{cor.Phi6} and~\ref{cor.Phi10}, and
a few more subtle
geometric properties illustrating~\iref{spec}{spec.geometry}
(the so called $Z$-splitting sections)
are
discusses in Subsection~\ref{s.Z-splitting}.

Finally, one has the following application to the fundamental
groups of plane curves, which were my original motivation for the
study of trigonal curves.

\theorem\label{th.plane}
Let $\DD\subset\Cp2$ be an irreducible plane curve with a singular
point of multiplicity $(\deg\DD-3)$. If $\DD$ admits a
$\GDG(\CQ)$-covering, then $\CQ$ is a quotient of one of the groups
$\CG3\oplus\CG3$, $\CG9$, $\CG5\oplus\CG5$, or $\CG7$.
\endtheorem

This theorem is proved in Subsection~\ref{proof.plane}. It is
worth mentioning that the fundamental group of any irreducible
plane curve~$\DD$ with a singular point of multiplicity
$(\deg\DD-1)$ is abelian, whereas for each
integer $m\ge1$ there is an
irreducible plane curve~$\DD$ with a singular point of
multiplicity $(\deg\DD-2)$ admitting a $\DG{4m+2}$-covering.

\subsection{Contents of the paper}\label{s.contents}
In Section~\ref{S.curves}, we
remind a few basic notions and facts related to trigonal curves in
Hirzebruch surfaces.
In Section~\ref{S.monodromy}, we
discuss the braid monodromy and Zariski--van Kampen theorem
computing the fundamental group, and then introduce the concept
and prove the existence of universal
trigonal curves related to a
prescribed monodromy group or a prescribed quotient
of~$\piaff{}$.
Section~\ref{S.DG} deals with uniform dihedral
quotients/coverings.
The principal result here is the proof of
Theorems~\ref{th.D} and~\ref{th.reducible}. In the course of the
proof, we treat the special case of
isotrivial curves, which are mostly
excluded from the consideration in the rest of the paper.
In Section~\ref{S.groups}, we discuss the geometric properties and
fundamental groups of the universal curves corresponding to
uniform dihedral coverings. The results are applied to illustrate
Statements~\iref{spec}{spec.geometry} and~\ditto{spec.larger}
and, in
particular, to prove Theorem~\ref{th.Oka}.
Finally, in Section~\ref{S.appl} we discuss some further
applications of the principal results: the relation to the
Mordell--Weil group, $Z$-splitting sections of trigonal curves,
and extensions to generalized trigonal curves and plane curves
with deep singularities.

\subsection{Acknowledgements}
I am grateful to A.~Klyachko for a number of valuable discussions
concerning the modular group, and to I.~Shimada, who brought to my
attention paper~\cite{Cox}.

\section{Trigonal curves\label{S.curves}}

We remind the basic notions and facts related to trigonal curves.
Important for the sequel are the notions of $m$-Nagata
equivalence, induced curves, and trigonal curves of torus type.

\subsection{Trigonal curves in Hirzebruch surfaces}\label{s.curves}
A \emph{Hirzebruch surface}~$\Sigma_d$
is a geometrically ruled rational surface with an
\emph{exceptional section}~$\EE$ of self-intersection $-d\le0$.
The \emph{fibers} of~$\Sigma_d$ are the fibers of the ruling
$\Sigma_d\to\Cp1$.
To avoid excessive notation,
we identify fibers and their images in the base~$\Cp1$.
The semigroup of classes of effective divisors on~$\Sigma_d$ is
freely generated by the classes $\ls|E|$ and~$\ls|F|$, where $F$
is any fiber.

A \emph{trigonal curve} is a reduced curve $\CC\subset\Sigma_d$
disjoint
from the exceptional section $\EE\subset\Sigma_d$ and
intersecting each fiber at three points; in other words,
$\CC\in\ls|3\EE+3dF|$.
A trigonal curve is called \emph{simple} if all its singular
points are simple, \ie, those of type~$\bA_p$, $p\ge1$, $\bD_q$,
$q\ge4$, $\bE_6$, $\bE_7$, or~$\bE_8$.

A \emph{singular fiber} of a trigonal curve $\CC\subset\Sigma_d$
is a fiber intersecting~$\CC$ geometrically at fewer than three
points.
For the topological types of singular fibers, we use the following
notation, referring to the types of the singular points of~$\CC$
(to simplify a few statements, we sometimes extend Arnol$'$d's
notation~$\bJ$, $\bE$ for the non-simple triple singular points,
see~\cite{AVG},
to the type $\bA$ and~$\bD$ points as well):
\def\Kodaira's{Kodaira's}
\def\Kodaira's{\ignorespaces}
\Dashes
\dash
$\tA{_0}=\tJ_{0,0}$ (\Kodaira's $\I_0$):
a nonsingular fiber;
\dash
$\tA{_0^*}=\tJ_{0,1}$ (\Kodaira's $\I_1$):
a simple vertical tangent;
\dash
$\tA{_0^{**}}=\tE_0$ (\Kodaira's $\II$):
a vertical inflection tangent;
\dash
$\tA{_1^*}=\tE_1$ (\Kodaira's $\III$):
a node of~$\CC$ with one of the branches vertical;
\dash
$\tA{_2^*}=\tE_2$ (\Kodaira's $\IV$):
a cusp of~$\CC$ with vertical tangent;
\dash
$\tA{_{p}}=\tJ_{0,p+1}$ (\Kodaira's $\I_{p+1}$), $p\ge1$,
\nodash
$\tD{_q}=\tJ_{1,q-4}$ (\Kodaira's $\I_{q-4}^*$), $q\ge4$,
\nodash
$\tE_6$ (\Kodaira's $\IV^*$), $\tE_7$ (\Kodaira's $\III^*$),
$\tE_8$ (\Kodaira's $\II^*$):
a simple singular point of~$\CC$ of the same type
with the minimal possible local intersection index with the fiber;
\dash
$\tJ{_{r,p}}$, $r\ge2$, $p\ge0$,
\nodash
$\tE{_{6r+\epsilon}}$, $r\ge2$, $\epsilon=0,1,2$:
a non-simple singular point of~$\CC$ of the same type.
\endDashes
For the simple singular fibers we
also list parenthetically Kodaira's notation for the corresponding
singular fiber of the covering elliptic surface. In the case of a
simple fiber, the
$\bA$--$\bD$--$\bE$ notation refers as well to the incidence
graph of $(-2)$-curves in the corresponding elliptic fiber;
this graph is an affine Dynkin diagram.

The fibers of type~$\tE$ (including $\tE_0$, $\tE_1$, and~$\tE_2$)
are called \emph{exceptional}.

\subsection{Nagata transformations}
A positive (negative) \emph{Nagata transformation} is the
birational transformation $\Sigma_d\dashrightarrow\Sigma_{d\pm1}$
consisting in blowing up a point~$P$ on (respectively, not on) the
exceptional section~$\EE$ and blowing down the proper transform of the
fiber through~$P$. An \emph{$m$-fold} Nagata transformation is a
sequence of $m$ Nagata transformations \emph{of the same sign}
over the same point of the base.

\definition\label{def.Nagata}
Two trigonal curves~$\CC$, $\CC'$ are called \emph{$m$-Nagata
equivalent} if $\CC'$ is the proper transform of~$\CC$
under a sequence of $m$-fold Nagata transformations. The
special case $m=1$ is referred to as just \emph{Nagata equivalence}.
\enddefinition

Pick a fiber~$F_0$ (\emph{fiber at infinity}),
consider the affine chart
$\Sigma_d\sminus(\EE\cup F_0)$, and choose affine coordinates
$(x,y)$ so that the fibers of the ruling be the vertical lines
$x=\const$. In these coordinates, any trigonal curve~$\CC$ is
given by an equation of the form
$\sum_{i=0}^3y^ib_i(x)$, where $b_i(x)$ is a polynomial of degree
up to $d(3-i)$.
In the same coordinates, a positive Nagata transformation is
given by the coordinate change
$$
x=x',\quad y=y'\!/x'\eqtag\label{eq.Nagata}
$$
(assuming that the image of the fiber contracted is the origin
$(x',y')=(0,0)$),
and the equation of the transform of~$\CC$ is obtained from the
original equation of~$\CC$
by the substitution and clearing the fractions. A negative Nagata
transformation is given by $x=x'$, $y=y'x'$. For the result to be
disjoint from~$\EE$, the original
curve~$\CC$ must have a triple singular point
at the origin (the blow-up center); then,
after the substitution, one can cancel $x^{\prime3}$.

Under a single positive Nagata transformation, the topological
type of the fiber affected changes as follows:
$$
\tJ_{r,p}\mapsto\tJ_{r+1,p},
 \ r,p\ge0,\qquad
\tE_{6r+\epsilon}\mapsto\tE_{6(r+1)+\epsilon},
 \ r\ge0,\ \epsilon=0,1,2.
\eqtag\label{eq.types}
$$
This statement can easily be obtained using~\eqref{eq.Nagata} and
the local normal forms.
In each series, there are exactly two simple
singularities, those with $r=0$ or~$1$. Each series starts with
its only type~$\tA$ singularity, corresponding to $r=0$.

\subsection{The $j$-invariant}\label{s.j}
The \emph{\rom(functional\rom) $j$-invariant}
$j_\CC\:\Cp1\to\Cp1$ of
a trigonal curve $\CC\subset\Sigma_d$ is defined as the analytic
continuation of the function sending a
nonsingular
fiber~$F$
to the $j$-invariant (divided
by~$12^3$)
of the elliptic
curve covering~$F$
and ramified at
$F\cap(\CC\cup\EE)$.
The curve~$\CC$
is called \emph{isotrivial} if $j_\CC=\const$.
Such curves are
easily enumerated, see Subsection~\ref{s.isotrivial}.

In appropriate affine coordinates $(x,y)$ as above
the curve~$\CC$ can
be given by its \emph{Weierstra{\ss} equation}
$$
y^3+3p(x)y+2q(x)=0.\eqtag\label{eq.W}
$$
Then, the $j$-invariant is given by
$$
j_\CC(x)=\frac{p^3}\Delta,\quad\text{where}\quad
 \Delta(x)=p^3+q^2.\eqtag\label{eq.j}
$$
Up to a constant factor, $\Delta(x)$ is the discriminant
of~\eqref{eq.W} with respect to~$y$.

By definition, $j_\CC$ is invariant under Nagata transformations.
Any holomorphic map $j\:\Cp1\to\Cp1$ is the $j$-invariant of a
certain trigonal curve~$\CC$, which is unique
up to isomorphism and Nagata equivalence.
Exceptional singular fibers of~$\CC$ are those where $j$ takes
value~$0$ or~$1$ and has ramification index $\ne0\bmod3$ or
$\ne0\bmod2$, respectively. Singular fibers of type~$\tJ_{r,0}$,
$r\ge1$, are not detected by the $j$-invariant, and all other
singular fibers of~$\CC$
(those of types~\smash{$\tJ_{r,p}$}, $r\ge0$, $p\ge1$)
are precisely those where $j$ takes value~$\infty$. At such a
fiber, the ramification index of~$j$ is~$p$.

Informally, the $j$-invariant~$j_\CC$ determines the
type~$\tJ_{r,p}$ or $\tE_{6r+\epsilon}$ of each singular fiber
of~$\CC$ up to a choice of the integer $r\ge0$. Thus, in order to
select a single curve in its Nagata equivalence class, one needs
to fix its \emph{type specification}, \ie, select a precise type
of each singular fiber and, possibly, assign types $\tJ_{r,0}$
(not detected by the $j$-invariant)
to a few generic fibers.

\subsection{Maximal curves and skeletons}\label{s.max}
A non-isotrivial trigonal curve~$\CC$
is called \emph{maximal} if it has the following properties:
\roster
\item\local{noD4}
$\CC$ has no singular fibers of type~$\tJ_{r,0}$, $r\ge1$;
\item\local{0,1,infty}
$j=j_{\CC}$ has no critical values other than~$0$, $1$, and~$\infty$;
\item\local{le3}
each point in the pull-back $j\1(0)$ has ramification index at
most~$3$;
\item\local{le2}
each point in the pull-back $j\1(1)$ has ramification index at
most~$2$.
\endroster
An important property of maximal trigonal
curves is their rigidity, see~\cite{degt.kplets}: any small
fiberwise equisingular deformation of such a curve
$\CC\subset\Sigma_d$ is
isomorphic to~$\CC$. Any maximal trigonal curve is defined over an
algebraic number field.

The $j$-invariant of a maximal trigonal curve~$\CC$ can be described
by its \emph{skeleton}, which is defined as the embedded bipartite
graph $\Sk_\CC:=j_\CC\1[0,1]\subset S^2\cong\Cp1$, with the
\black-- and \white-vertices being the pull-backs of~$0$ and~$1$,
respectively. By definition, all \black-- (respectively, \white--)
vertices of~$\Sk_\CC$ are of valency~$\le3$
(respectively,~$\le2$); in the drawings, we omit bivalent
\white-vertices, assuming such a vertex at the center of each edge
connecting two \black-vertices. Each connected component of the
complement $S^2\sminus\Sk_\CC$ is a topological disk; hence,
instead of the embedding $\Sk_\CC\subset S^2$ one can
regard~$\Sk_\CC$ as a bipartite ribbon graph of genus zero.

Each skeleton $\Sk$ as above gives rise to a topological ramified
covering $S^2\to\Cp1$. By the Riemann existence theorem, the
latter is realized by a holomorphic map $\Cp1\to\Cp1$, unique up
to a M\"obius transformation in the source. Hence, the skeleton
$\Sk_\CC$ determines $j_\CC$. The type specification of a maximal
trigonal curve~$\CC$
can be regarded as a function assigning an integer $r\ge0$
to each \black-vertex of valency~$\le2$, each \white-vertex of
valency~$1$, and each region of~$\Sk_\CC$.

\subsection{Induced curves}
Consider a Hirzebruch surface $\Sigma:=\Sigma_d$ over a base
$\BB\cong\Cp1$,
let $\BB'\cong\Cp1$ be another rational curve,
and let $\tj\:\BB'\to\BB$ be a nonconstant
holomorphic map.
Then the ruled surface
$\Sigma':=\tj^*\Sigma$ over~$\BB'$ is also a Hirzebruch surface;
it is isomorphic to $\Sigma_{d\cdot\deg\tj}$. Furthermore, given a
trigonal curve $\CC\subset\Sigma$, its divisorial pull-back
$\CC':=\tj^*\CC\subset\Sigma'$ is also a trigonal curve; it is
said to be \emph{induced} from~$\CC$ by~$\tj$ or obtained
from~$\CC$ by a \emph{rational base change}.

In appropriate affine coordinates $(x,y)$ in~$\Sigma$ and
$(x',y')$ in~$\Sigma'$, the ramified covering $\Sigma'\to\Sigma$
is the map
$$
(x',y')\mapsto(x,y)=
 \left(\frac{u(x')}{v(x')},\frac{y'}{v^d(x')}\right);
 \eqtag\label{eq.change}
$$
here, $x$ and~$x'$ are affine parameters in~$\BB$ and~$\BB'$,
respectively, and $\tj$ is given by the reduced fraction
$\tj(x')=u(x')/v(x')$. An equation for~$\CC'$ is obtained from
that for~$\CC$ by substituting~\eqref{eq.change} and clearing
denominators.

Locally, the substitution is given by
$(x',y')\mapsto(x,y)=(x^{\prime m},y')$, where $m$ is the
ramification index of~$\tj$ at $x'=0$, and, using local normal
forms, one can easily find the types of the singular fibers
of~$\CC'$ in terms of those of~$\CC$. Next lemma characterizes
Kodaira type~$\I$ fibers, \ie, types $\tA_0^*$ and $\tA_p$,
$p\ge1$.

\lemma\label{type.I}
If a trigonal curve~$\CC$ has Kodaira type~$\I$ singular fibers
only, then any curve~$\CC'$ induced from~$\CC$ by a rational base
change~$\tj$ is simple\rom; in fact, $\CC'$ also has Kodaira type~$\I$
singular fibers only.
\endlemma

\proof
The proof is a simple computation using local normal forms, as
explained above. Assume that $\tj$ sends $F'=\{x'=0\}$ to
$F=\{x=0\}$ and has ramification index~$m$ at~$F'$. If $F$ is of
type~$\I_{p+1}$, $p\ge0$, then $F'$ is of type $\I_{m(p+1)}$. If
$F$ is of any other type, its local normal form is $y^3+xa(x,y)$
and, for any $m>6$, the induced fiber~$F'$ is not simple.
\endproof

\subsection{Trigonal curves of torus type}\label{s.torus}
A trigonal curve $\CC$ in an even Hirzebruch surface~$\Sigma_{2k}$
is said to be of \emph{torus type} if there are
sections $f_i$ of $\CO(\EE+ikF)$, $i=0,2,3$,
such that $\CC$ is the zero set of the section
$f_2^3+f_0f_3^2$. Informally, in affine coordinates $(x,y)$ with
$\EE=\{y=\infty\}$, the equation of~$\CC$ has the form
$$
[y+a_2(x)]^3+[a_1(x)y+a_3(x)]^2=0\eqtag\label{eq.torus}
$$
for some polynomials $a_i(x)$ of degree up to $ki$, $i=1,2,3$.
Each representation of the equation of~$\CC$ in this form
(up to an obvious equivalence)
is called a
\emph{torus structure} on~$\CC$.

\lemma\label{torus.invariant}
Torus structures are
preserved under rational
base changes and $2$-fold Nagata transformations.
\endlemma

\proof
Consider a rational base change $\Cp1\to\Cp1$ given by
$x=\tj(x'):=u(x')/v(x')$, $y=y'/v^{2k}(x')$,
see~\eqref{eq.change}.
Substituting to~\eqref{eq.torus} and
clearing denominators, one obtains
$$
[y'+a_2'(x')]^3+[a_1'(x')y'+a_3'(x')]^2=0,\eqtag\label{eq.torus'}
$$
where
$a_i'(x')=a_i(u/v)v^{ki}(x')$, $i=1,2,3$, \ie, again an equation
of the form~\eqref{eq.torus}.

Now, consider a positive $2$-fold Nagata transformation $x=x'$,
$y=y'/x^{\prime2}$, see~\eqref{eq.Nagata};
we can assume that the fiber
contracted is over $x=0$. Substituting and clearing denominators,
one obtains~\eqref{eq.torus'} with $a_i'(x')=a_i(x')x^{\prime i}$,
$i=1,2,3$.

Finally, consider a negative $2$-fold Nagata transformation; in
appropriate affine coordinates it is given by
$x=x'$, $y=y'x^{\prime2}$. For the
transform to be disjoint from the exceptional section,
the singularity of the original curve at the origin must be
adjacent to $\bJ_{2,0}$, \ie, all terms
$x^\Ga y^\Gb$, $2\Ga+\Gb<6$, of the original equation must vanish.
Evaluating at $x=0$, one can easily see that
$a_1(0)=a_2(0)=a_3(0)=0$, and then, step by step, one can conclude
that $a_i(x)=x^ia_i'(x)$ for some polynomials $a_i'(x)$,
$i=1,2,3$. Substituting and cancelling $x^{\prime6}$, one arrives
at~\eqref{eq.torus'}.
\endproof

\Remark\label{rem.torus}
If $\CC$ is the proper transform of a curve~$\CC'$
of torus type under a negative Nagata
transformation contracting a fiber~$F$, then the
\emph{divisorial} transform
of~$\CC'$ inherits a torus
structure: one has $\CC+sF=\{f_2^3+f_0f_3^2=0\}$
for some $s\ge0$. What is shown in the proof above is that, if
after a $2$-fold transformation $\CC$ is still disjoint
from~$\EE$, then $s=6$ and the curves $\{f_2=0\}$ and
$\{f_3=0\}$ contain~$F$ with multiplicity at least~$2$ and~$3$,
respectively, so that $6F$ can be factored out.
\endRemark

\lemma\label{torus.monodromy}
If a trigonal curve $\CC\subset\Sigma_{2k}$ is of torus type,
there exists a triple covering $X\to\Sigma_{2k}$ ramified at
$\CC\cup\EE$ with the full monodromy group~$\SG3$.
\endlemma

\proof
In affine coordinates $(x,y,z)$ in $\C^2\times\C$, the covering
surface~$X$ is given by
$z^3+3(y+a_2)+2(a_1y+a_3)=0$, see~\eqref{eq.torus}. More formally,
$X$ is the normalization of the triple section of $\CO(\EE+kF)$
given by $z^3+3f_0f_2z+2f_0^2f_3=0$. Restricting to a generic
fiber $x=\const$, one obtains a covering $\CC_x\to\Cp1$, where
$\CC_x\subset\Sigma_1$ is a trigonal curve with the set of
singular fibers $\tA_1^*\oplus3\tA_0^*$. All such curves are
deformation equivalent (essentially, they are nodal plane cubics
projected from a generic point in the tangent to one of the
branches at the node), and their monodromy groups are easily shown
to equal~$\SG3$, see Remark~\ref{rem.S3.torus} below for details.
\endproof

\Remark
If $\CC$ is irreducible, the converse of
Lemma~\ref{torus.monodromy} also holds: it is given by
Theorem~\ref{th.Oka}, proved in Subsection~\ref{proof.Oka}.
\endRemark

\section{The braid monodromy\label{S.monodromy}}

In this section, we define the braid monodromy and its various
reductions, cite Zariski--van Kampen theorem and its implications
for the particular case of trigonal curves, introduce the concept
of universal curves, and prove their existence.

\subsection{Groups to be considered}\label{s.groups}
Let
$\CH=\Z a\oplus\Z b$ be a rank~$2$ free abelian group
with the skew-symmetric bilinear form
$\bigwedge^2\CH\to\Z$ given by $a\cdot b=1$. We fix the
notation~$\CH$, $a$, $b$ throughout the paper
and \emph{define}
$\tMG:=\SL(2,\Z)$
as the group $\Sp\CH$
of symplectic auto-isometries of~$\CH$; it is
generated by the
isometries $\X,\Y\:\CH\to\CH$ given (in the basis $\{a,b\}$ above)
by the matrices
$$
\X=\bmatrix-1&1\\-1&0\endbmatrix,\qquad
\Y=\bmatrix0&-1\\1&\phantom{-}0\endbmatrix.
$$
One has $\X^3=\id$ and $\Y^2=-\id$.
The \emph{modular group} $\MG:=\PSL(2,\Z)$ is
the
quotient $\tMG/\!\pm\id$. We retain the notation $\X$, $\Y$ for
the generators of~$\MG$. One has
$$
\MG=\<\X\,|\,\X^3=1\>\mathbin*\<\Y\,|\,\Y^2=1\>\cong\CG3\mathbin*\CG2.
$$

The \emph{braid group} $\BG3$ is the group
$$
\BG3=\<\Gs_1,\Gs_2\,|\,\Gs_1\Gs_2\Gs_1=\Gs_2\Gs_1\Gs_2\>=
 \<u,v\,|\,u^3=v^2\>,
$$
where $u=\Gs_2\Gs_1$ and $v=\Gs_2\Gs_1^2$. The center $Z(\BG3)$ is
the infinite cyclic subgroup generated by $u^3=v^2$, and the quotient
$\BG3/Z(\BG3)$ is isomorphic to~$\MG$. We define the epimorphism
$\BG3\onto\tMG$ (and further to~$\MG$) \via
$$
\Gs_1\mapsto\X\Y=\bmatrix1&1\\0&1\endbmatrix,\quad
 \Gs_2\mapsto\X^2\Y\X\1=\bmatrix\phantom{-}0&1\\-1&2\endbmatrix.
\eqtag\label{eq.BG3toMG}
$$
Then $u\mapsto-\X\1$ and $v\mapsto-\Y$.
This unusual
choice of generators is explained in Remark~\ref{modular.bm} below.

The braid group~$\BG3$ acts on the free group
$\<\Ga_1,\Ga_2,\Ga_3\>$ \via
$$
\Gs_1\:\Ga_1\mapsto\Ga_1\Ga_2\Ga_1\1,\quad
 \Ga_2\mapsto\Ga_1;\qquad
\Gs_2\:\Ga_2\mapsto\Ga_2\Ga_3\Ga_2\1,\quad
 \Ga_3\mapsto\Ga_2.
$$
According to E.~Artin~\cite{Artin}, $\BG3$ can be identified with
the group of automorphisms of $\<\Ga_1,\Ga_2,\Ga_3\>$ taking each
generator to a conjugate of a generator and preserving the product
$\Gr:=\Ga_1\Ga_2\Ga_3$. In what follows, we fix the notation~$\AA$
for the free group $\<\Ga_1,\Ga_2,\Ga_3\>$ supplied with  this
$\BG3$-action. We do not distinguish between the original basis
$\{\Ga_1,\Ga_2,\Ga_3\}$ and any other basis in its $\BG3$-orbit;
any basis in the latter orbit is called \emph{geometric}.
We also reserve the notation~$\Gr$ for
the product $\Ga_1\Ga_2\Ga_3$.

The group~$\AA$ is also equipped with a distinguished homomorphism
$\deg\:\AA\to\Z$, $\Ga_i\mapsto1$, which does not depend on the
choice of a geometric basis.


Any finite index subgroup $G\subset\BG3$ must intersect the center
$Z(\BG3)$. We define the
\emph{depth} $\depth G$ as the minimal positive integer~$m$ such
that $(\Gs_2\Gs_1)^{3m}\in G$.

We are mainly interested in subgroups of~$\BG3$, $\tMG$, or~$\MG$
up to conjugation. Given two subgroups $H,H'\subset G$, we say
that $H'$ is \emph{subconjugate} to~$H$,
notation $H'\prec H$,
if $H'$ is conjugate in~$G$ to a subgroup of~$H$.

\subsection{Skeletons and genus}\label{s.Sk}
Given a finite index subgroup $G\subset\MG$, the quotient
$\MG\!/G$ can be given a natural structure of a bipartite
ribbon graph; it is denoted by $\Sk_G$ and called the
\emph{skeleton} of~$G$. The set of edges of~$\Sk_G$ is $\MG\!/G$.
This set has a canonical left $\MG$-action, and we define the
\black-- and \white-vertices of~$\Sk_G$ as the orbits of~$\X$
and~$\Y$, respectively. The cyclic order (the ribbon graph
structure) at a trivalent \black-vertex is given by~$\X\1$; all
other vertices are at most bivalent, and the cyclic order is
unique.
Alternatively, there is a natural $\MG$-action on the
infinite \emph{Farey tree}~$\mathrm{F}$
(the only bipartite tree with all
\black-vertices of valency~$3$ and all \white-vertices of
valency~$2$), and one can define $\Sk_G$ as $\mathrm{F}\!/G$.
For more details and all proofs, see~\cite{monodromy}.

By definition,
the valency of each \black-- (respectively, \white--) vertex
of~$\Sk_G$ is divisible by~$3$ (respectively,~$2$). In the
drawings, we omit bivalent \white-vertices, assuming that such a
vertex is to be inserted at the middle of each edge connecting two
\black-vertices.
Conversely, the set of edges of any bipartite
ribbon graph~$\Sk$ satisfying the above valency condition admits a
natural $\MG$-action, and the stabilizer of a fixed edge of~$\Sk$
is a certain finite index subgroup $G\subset\MG$.
One has $\Sk\cong\Sk_G$.

\Remark
Both skeletons of maximal trigonal curves, see
Subsection~\ref{s.max}, and skeletons of subgroups of~$\MG$ defined
above are bipartite ribbon graph with all \black-vertices of
valency~$\le3$ and all \white-vertices of valency~$\le2$;
in the figures, we use the same convention about bivalent
\white-vertices.
Note though that there is
a certain difference between the two classes: skeletons of
trigonal curves are required to be of genus zero (due to the fact
that we consider curves over a rational base only), whereas
skeletons of subgroups are \emph{not} allowed to have bivalent
\black-vertices.
Still, the two notions are closely related, \cf.
Remark~\ref{rem.group->curve} below.
\endRemark

The monovalent vertices of~$\Sk_G$ are in a one-to-one
correspondence with the torsion elements of~$G$, and its
\emph{regions}, \ie, minimal left turn cycles, correspond to the
conjugacy classes of indivisible parabolic elements of~$G$.
If $G$ is torsion free,
equivalently,
if $\Sk_G$ has no monovalent vertices,
there is a canonical isomorphism $G=\pi_1(\Sk_G)$. In general, $G$
is isomorphic to an appropriately defined
orbifold fundamental group
$\pi^{\mathrm{orb}}_1(\Sk_G)$.

\definition\label{def.genus}
If $[\MG:G]<\infty$, the genus of the minimal surface
supporting the skeleton~$\Sk_G$
is called the \emph{genus} of~$G$. The genus of
a finite index subgroup of~$\tMG$ or~$\BG3$ is defined as the
genus of its image in~$\MG$.
\enddefinition

Genus is nonnegative,
invariant under conjugation, and
monotonous: one has
$\operatorname{genus}(G')\ge\operatorname{genus}(G)\ge0$
whenever $G'\prec G$.

Definition~\ref{def.genus} is equivalent to the classical
definition of genus, see Remark~\ref{rem.genus} below. It is worth
emphasizing that, when speaking about a subgroup $G\subset\BG3$ of
genus zero, we always assume that $G$ is of finite index, \ie,
$G\cap Z(\BG3)\ne\{1\}$.

\subsection{Proper sections and braid monodromy}\label{s.sections}
In the exposition below, we follow the approach of~\cite{degt.e6},
which makes certain choices in the definition of braid monodromy
more canonical. We refer to~\cite{degt.e6} for most proofs, which
are omitted.

Fix a Hirzebruch surface~$\Sigma_d$
and a trigonal curve $\CC\subset\Sigma_d$.
The term `section' stands for a continuous section of (a
restriction of) the fibration $p\:\Sigma_d\to\Cp1$.
For any fiber~$F$ of~$\Sigma_d$, the complement
$F^\circ:=F\sminus\EE$ is an affine space over~$\C$. Hence, one
can speak about the convex hull of a subset of~$F^\circ$. For a
subset $S\subset\Sigma_d\sminus\EE$, we denote by $\Conv_F S$
the convex hull of $S\cap F\scirc$ in~$F\scirc$ and
let $\Conv S$ be the union of $\Conv_FS$ over all fibers~$F$.

\definition\label{def.proper}
Let $\disk\subset\Cp1$ be a closed topological disk. A
section $s\:\disk\to\Sigma_k$ of~$p$ is called \emph{proper} if its
image is disjoint from both~$\EE$ and $\Conv\CC$.
\enddefinition

Any disk $\disk\subset\Cp1$ admits a proper section
$s\:\disk\to\Sigma_k$, unique up to homotopy in the class
of proper sections.

Fix a disk $\disk\subset\Cp1$ and
let $F_1,\ldots,F_r\in\disk$ be all singular and, possibly, some
nonsingular fibers of~$\CC$ that belong to~$\disk$.
Assume that
all these fibers are in the interior of~$\disk$.
Let
$\disk\scirc=\disk\sminus\{F_1,\ldots,F_r\}$ and
fix a \emph{reference fiber}
$F\in\disk\scirc$. Then, given a proper
section~$s$, one can define the group
$\pi_F:=\pi_1(F\scirc\sminus\CC,s(F))$ and the
\emph{braid monodromy}, which is the anti-homomorphism
$\bm\:\pi_1(\disk\scirc,F)\to\Aut\pi_F$ sending a loop~$\Gg$ to
the automorphism obtained by dragging~$F$ along~$\Gg$ and keeping
the reference point in~$s$.

\definition
Let~$D$ be an oriented punctured disk, and let $b\in\partial D$. A
\emph{geometric basis} in~$D$ is a basis $\{\Gg_1,\ldots,\Gg_r\}$
for the free group $\pi_1(D,b)$ formed by the classes of
positively oriented lassoes about the punctures,
pairwise disjoint except
at the common reference point~$b$ and such that
$\Gg_1\ldots\Gg_r=[\partial D]$.
\enddefinition

Shrink the reference fiber~$F$ to a closed disk containing
$\Conv_F\CC$ in its interior and $s(F)$ in its boundary. Pick a
geometric basis for~$\pi_F$ and identify it with a geometric basis
$\{\Ga_1,\Ga_2,\Ga_3\}$
for~$\AA$, establishing an isomorphism $\pi_F\cong\AA$.
Due to Artin's theorem~\cite{Artin},
under this
isomorphism the braid monodromy~$\bm$ takes values in the braid
group $\BG3\subset\Aut\AA$, which explains the term. The
$\BG3$-valued braid monodromy~$\bm$ thus defined is independent of
the choice of a proper section, and another choice
of the geometric bases
used for the identification $\pi_F=\AA$ results in the
global conjugation by a fixed braid $\Gb\in\BG3$, \ie, in the map
$\Gg\mapsto\Gb\1\bm(\Gg)\Gb$.

\lemma\label{infty}
Assume that the disk~$\disk$ contains all singular fibers
of~$\CC$. Then, in any geometric basis $\{\Gg_1,\ldots,\Gg_r\}$ in
$\disk\scirc$,
the so called \emph{monodromy at infinity}
$\bm(\Gg_1\ldots\Gg_r)=\bm[\partial\disk]$
equals $(\Gs_2\Gs_1)^{3d}$\rom;
it is the conjugation by~$\Gr^d$.
\qed
\endlemma

In what follows, we will take for
the disk~$\disk$ the complement of a
small regular neighborhood of a nonsingular fiber $F_0\in\Cp1$.
Due to Lemma~\ref{infty}, the braid monodromy over~$\disk$
factors to an anti-homomorphism
$\bm\:\pi_1(\BB\scirc)\to\BG3/(\Gs_2\Gs_1)^{3d}$, where
$\BB\scirc=\Cp1\sminus\bigcup_{i=1}^rF_r$. This map is independent
of the choice of~$F_0$. Its image $\BM_\CC(\BG3)$, regarded as a
subgroup of~$\BG3$, is called the \emph{monodromy group} of~$\CC$;
it is defined by~$\CC$ up to conjugation. We will also consider
the reductions of~$\bm$ to the groups~$\tMG$ and~$\MG$;
their images are denoted by
$\BM_\CC(\tMG)$ and $\BM_\CC(\MG)$, respectively.

\Remark\label{modular.bm}
The $\MG$-valued braid monodromy $\bm\:\pi_1(\BB\scirc)\to\MG$ is
always well defined. For a non-isotrivial
trigonal curve~$\CC$, it can be
expressed in terms of its $j$-invariant~$j_\CC$.
Let $\BB_\MG$ be
the Riemann sphere $\C\cup\{\infty\}\cong\Cp1$, and denote
$\BB_\MG\scirc=\BB_\MG\sminus\{0,1,\infty\}$. The
identification of $\BB_\MG$
with the \emph{modular curve}
$\MG\backslash\fH^*$, see, \eg,~\cite{Shimura},
gives rise to a canonical principal
$\MG$-bundle
over~$\BB_\MG\scirc$, which defines the
\emph{monodromy \hbox{\rom(anti-\rom)}\penalty0representation}
$\bm_\MG\:\BB_\MG\to\MG$. Then, up to global conjugation, the
$\MG$-valued braid monodromy of a curve~$\CC$ is the composition
$\bm_\MG\circ(j_\CC)_*$.
This fact is well known for the homological invariant of an
elliptic surface;
the relation between the modular representation
and the braid monodromy of a trigonal curve can be established,
\eg, using the computation in~\cite{degt.kplets}.

Accidentally, it is this reduction that motivates our not quite
usual choice of the epimorphism $\BG3\onto\tMG$ fixed in
Subsection~\ref{s.groups}: it is consistent with a \emph{canonical
basis} $\{\Ga_1,\Ga_2,\Ga_3\}$ for~$\pi_F$ constructed
in~\cite{degt.kplets} and the basis
$a=\Ga_2\Ga_1$, $b=\Ga_1\Ga_3$ for the group~$\CH$ regarded as the
$1$-homology of the double covering of~$F$, see
Subsection~\ref{s.quotients} and~\eqref{eq.CH} below.
\endRemark

\Remark\label{rem.group->curve}
It follows, in particular, that the skeleton of a maximal trigonal
curve without type $\tE_{6r+2}$, $r\ge0$,
singular fibers, \ie, without
bivalent \black-vertices, can be identified with the skeleton of
its monodromy group, see~\cite{monodromy}. Hence, such a curve is
determined by the conjugacy class of its $\MG$-valued monodromy
group up to Nagata equivalence.
Furthermore, each genus zero subgroup of~$\MG$ determines a unique
Nagata equivalence class of maximal trigonal curves; the curve can
be made unique up to isomorphism if all singular fibers are
required to be of type~$\tA$.
\endRemark

\lemma\label{monodromy.Nagata}
In the notation of Lemma~\ref{infty}, assume that $\CC'$ is
obtained from~$\CC$ by a positive Nagata transformation at a
fiber~$F_i$, and let $\bm'$ be the braid monodromy of~$\CC'$. Then
$\bm'(\Gg_j)=\bm(\Gg_j)$ for $j\ne i$, and
$\bm'(\Gg_i)=\bm(\Gg_i)(\Gs_2\Gs_1)^3$.
\endlemma

\proof
The statement follows from~\eqref{eq.Nagata}: the
monodromy of~$\CC'$ about~$F_i$ differs from that of~$\CC$ by an
extra full twist, \ie, $(\Gs_2\Gs_1)^3$ or the conjugation by~$\Gr$.
\endproof

\Remark\label{rem.type}
According to Lemma~\ref{monodromy.Nagata}, the choice of a type
specification selecting a curve within its Nagata equivalence
class is equivalent to the choice of a lift to~$\BG3$ of the
$\MG$-valued braid monodromy determined by the $j$-invariant.
\endRemark

\lemma\label{monodromy.induced}
Assume that a trigonal curve~$\CC'$ is induced from~$\CC$ by a
rational base change $\tj\:\BB'\to\BB$. Then
$\BM_{\CC'}(\BG3)\prec\BM_\CC(\BG3)$.
\endlemma

\proof
The data $\Delta'$, $s'$, $F'$ necessary to define the
monodromy~$\bm'$ of~$\CC'$ can be pulled back from those
for~$\CC$; namely, let $\Delta'=\tj\1(\Delta)$ (formally, with a
few extra cuts to make it a disk) and $s'=\tj^*s$ and take
for~$F'$ any fiber in the pull-back of~$F$.
Then obviously $\bm'=\bm\circ\tj_*$, and the statement follows.
\endproof

\subsection{The Zariski--van Kampen theorem}\label{s.vanKampen}
Consider a trigonal curve $\CC\subset\Sigma_d$.
Keeping the notation of Subsection \ref{s.sections}, fix a
nonsingular fiber~$F_0$ of~$\CC$, denote by $F_1$,\dots, $F_r$ all
its singular
(and possibly some nonsingular) fibers,
identify the group of the
reference fiber with~$\AA$, and let $\bm\:\pi_1(\disk\scirc)\to\BG3$
be the resulting braid monodromy.
Fix also a geometric basis $\{\Gg_1,\ldots,\Gg_r\}$ in the
punctured disk~$\disk\scirc$ and let $\Gg_0=(\Gg_1\ldots\Gg_r)\1$:
this class is represented by a small loop about~$F_0$.

The following theorem is essentially contained
in~\cite{vanKampen}.

\theorem\label{th.vanKampen}
In the notation above, one has
$$
\multline\textstyle
\pi_1(\Sigma_d\sminus(\CC\cup\EE\cup\bigcup_{i=0}^rF_r))=
\bigl\<\Ga_1,\Ga_2,\Ga_3,\tGg_0,\tGg_1,\ldots,\tGg_r\bigm|\\
 \tGg_i\1\Ga_j\tGg_i=\bm_i(\Ga_j),\ i=1,\ldots,r,\ j=1,2,3,\
 \tGg_0\tGg_1\ldots\tGg_r\Gr^d=1\bigr\>,
\endmultline
$$
where $\bm_i=\bm(\Gg_i)$, $i=1,\ldots,r$, and $\tGg_i$ is a certain
lift of $\Gg_i$, $i=0,\ldots,r$.
\qed
\endtheorem

Patching back in a fiber~$F_i$, $i=0,\ldots,r$, results in an
extra relation $\tGg_i=1$. Hence, one has the following corollary.

\corollary\label{groups}
One has
$$
\align
\piaff\CC&=
\bigl\<\Ga_1,\Ga_2,\Ga_3\bigm|
 \bm(\Gg_i)=\id,\ i=1,\ldots,r\bigr\>,\\
\piproj\CC&=
\bigl\<\Ga_1,\Ga_2,\Ga_3\bigm|
 \bm(\Gg_i)=\id,\ i=1,\ldots,r,\
 \Gr^d=1\bigr\>,
\endalign
$$
where each \emph{braid relation} $\bm(\Gg_i)=\id$ should be
regarded as a triple of relations
$\bm(\Gg_i)(\Ga_j)=\Ga_j$, $j=1,2,3$, or, equivalently, as
a set of relations
$\bm(\Gg_i)(\Ga)=\Ga$ for each $\Ga\in\<\Ga_1,\Ga_2,\Ga_3\>$.
\qed
\endcorollary

\definition
The presentations of $\piaff\CC$ and $\piproj\CC$
given by Corollary~\ref{groups} are called
\emph{geometric}. More precisely, a geometric presentation is an
epimorphism $\AA\onto\piaff\CC$ or $\AA\onto\piproj\CC$ obtained
by identifying a geometric basis in a nonsingular fiber~$F_0$ and
a geometric basis of~$\AA$.
\enddefinition

\Remark\label{rem.geometric}
In other words, Corollary~\ref{groups} states that
$\piaff\CC$ is the quotient of~$\AA$
by the minimal normal subgroup
containing
$g(\Ga)\Ga\1$ for all $g\in\BM_\CC(\BG3)$
and $\Ga\in\AA$, and $\piproj\CC$ is the further quotient of
$\piaff\CC$ by the
normal subgroup generated by $\Gr^d\in\piaff\CC$.
\endRemark

\corollary\label{commutants}
Given a trigonal curve $\CC\subset\Sigma_d$,
the canonical epimorphism
$\piaff\CC\onto\piproj\CC$ is a central extension
by
the infinite cyclic subgroup generated by $\Gr^d\in\piaff\CC$.
This epimorphism
induces an isomorphism
$[\piaff\CC,\piaff\CC]=[\piproj\CC,\piproj\CC]$.
In particular,
the group
$\piaff\CC$ is abelian if and only
if so is $\piproj\CC$.
\endcorollary

\proof
The element~$\Gr^d$ that
normally generates the kernel, see Remark~\ref{rem.geometric},
is central due to Lemma~\ref{infty}.
It follows from Corollary~\ref{groups} that
$\Gr^d$ is
mapped to an infinite order element of the abelianization
$\piaff\CC/[\piaff\CC,\piaff\CC]$. (In particular, $\Gr^d$ itself
is of infinite order and the kernel $\<\Gr^d\>$
is infinite.)
Hence, $\<\Gr^d\>\cap[\piaff\CC,\piaff\CC]=\{1\}$,
and the induced homomorphism of the commutants
is one-to-one.
Since it is also onto (as induced by an epimorphism), it is an
isomorphism.
\endproof

\corollary\label{cor.Nagata}
If two trigonal curves~$\CC$
and~$\CC'$ are $m$-Nagata equivalent, the quotients of the groups
$\piaff\CC$ and $\piaff{\CC'}$ by the commutator subgroups
$[\piaff\bullet,\<\Gr^m\>]$ are canonically
isomorphic.
\endcorollary

\proof
Computing the groups modulo the extra relations $[\Ga,\Gr^m]=1$,
$\Ga\in\AA$, one can
reduce the braid monodromy to $\BG3/(\Gs_2\Gs_1)^{3m}$. According
to
Lemma~\ref{monodromy.Nagata}, the two reductions coincide.
\endproof

\lemma\label{cor.induced}
Assume that a trigonal curve $\CC'$ is induced from
a curve~$\CC$ by a
rational base change $\tj\:\BB'\to\BB$. Then $\tj$ induces
epimorphisms $\piaff{\CC'}\onto\piaff\CC$ and
$\piproj{\CC'}\onto\piproj\CC$ compatible with geometric
presentations of the groups.
\endlemma

\proof
The statement
follows from
Corollary~\ref{groups} (see also
Remark~\ref{rem.geometric}) and Lemma~\ref{monodromy.induced}.
\endproof

\subsection{Universal curves}\label{s.universal}
Recall that there is a canonical faithful discrete action of the
modular group~$\MG$ on the upper half plane
$\fH:=\{z\in\C\,|\,\Im z>0\}$, so that $\MG\backslash\fH=\C$. Let
$\fH\scirc$ be $\fH$ with the orbits of~$i$ and $\exp(2\pi i/3)$
removed; the action is free on~$\fH\scirc$. For a subgroup
$G\subset\MG$, denote $\BB_G\scirc=G\backslash\fH\scirc$, and let
$j_G\:\BB_G\scirc\to\BB_\MG\scirc=\C\sminus\{0,1\}$ be the
projection. One obviously has $\deg j_G=[\MG:G]$ and
the image
$(j_C)_*\pi_1(\BB_G\scirc)\subset\pi_1(\BB_\MG\scirc)$ is the
pull-back of~$G$ under the modular representation, see
Remark~\ref{modular.bm}.
If $G$ is of finite index,
the Riemann surface~$\BB_G\scirc$
can be compactified to a ramified covering
$j_G\:\BB_G\to\BB_\MG=\Cp1$. In this case, one has
$\BB_G=G\backslash\fH^*$, where $\fH^*=\fH\cup\Q\cup\{\infty\}$,
see, \eg, \cite{Shimura}.

\Remark\label{rem.genus}
The skeleton~$\Sk_G$ of a finite index subgroup $G\subset\MG$ can
be identified with the ribbon graph $j_G\1[0,1]\subset\BB_G$, the
\black-- and \white-vertices being the pull-backs of~$0$ and~$1$,
respectively. Furthermore, since $j_G$ is ramified over~$0$, $1$,
and~$\infty$ only,
$\BB_G$ is a minimal surface
supporting~$\Sk_G$ and the genus of~$G$ given by
Definition~\ref{def.genus} equals the genus of~$\BB_G$, which is
the classical definition.
\endRemark

\lemma\label{j.factor}
Let $\CC$ be a non-isotrivial
trigonal curve over~$\BB$,
and assume that
$\BM_\CC(\MG)\prec G$
for some subgroup $G\subset\MG$.
Then the $j$-invariant
$j_\CC\:\BB\to\Cp1$ factors through $j_G\:\BB_G\to\BB_\MG=\Cp1$.
\endlemma

\proof
Assume that
$\BM_\CC(\MG)\subset G$ and
consider the (ramified) coverings $j_\CC\:\BB\to\Cp1$ and
$j_G\:\BB_G\scirc\to\BB_\MG\scirc=\Cp1\sminus\{0,1,\infty\}$.
Denote by $\BB_\MG\bcirc$ the base
$\BB_\MG\scirc$ with the critical values
of~$j_\CC$ removed, and let
$\BB_G\bcirc=j_G\1(\BB_\MG\bcirc)$ and
$\BB\bcirc=j_\CC\1(\BB_\MG\bcirc)$, so that the restrictions
of both~$j_G$ and~$j_\CC$ are unramified.
Under an appropriate choice of the base points,
the image
$\bm_\MG\circ(j_\CC)_*(\pi_1(\BB\bcirc))=\BM_\CC(\MG)\subset\MG$,
see Remark~\ref{modular.bm},
is a subgroup of~$G$.
Hence
$(j_\CC)_*(\pi_1(\BB\bcirc))\subset(j_G)_*(\pi_1(\BB_G\bcirc))=\bm_\MG\1(G)$
and
there is a lift $\tj\:\BB\bcirc\to\BB_G\bcirc$
splitting~$j_\CC$. In particular,
$[\MG:G]=\deg j_G\le\deg j_\CC<\infty$
and the Riemann surface~$\BB_G$ is well defined;
compactifying all curves,
one obtains a
desired
splitting $\BB_C\to\BB_G\to\BB_\MG=\Cp1$.
\endproof

\corollary\label{genus.0}
The monodromy groups
$\BM_\CC(\MG)$,
$\BM_\CC(\tMG)$, and $\BM_\CC(\BG3)$
of a non-isotrivial trigonal curve $\CC\subset\Sigma_d$
are subgroups of genus zero. In particular, they are of finite
index.
\endcorollary

\proof
It suffices to let $G=\BM_\CC(\MG)$ in Lemma~\ref{j.factor}.
The subgroup $\BM_\CC(\BG3)\subset\BG3$ is also
of finite index as $\BM_\CC(\BG3)\cap Z(\BG3)\ne\{1\}$,
see Lemma~\ref{infty}.
\endproof

\Remark
Note that the monodromy group $\BM_\CC(\MG)$ of an isotrivial
trigonal curve~$\CC$ is always finite cyclic,
hence of infinite index, see Subsection~\ref{s.isotrivial}.
\endRemark

\definition\label{def.universal}
A curve $\CC_G\subset\Sigma_d\to\BB_G$
is called a
\emph{universal trigonal curve} corresponding to the subgroup
$G\subset\MG$ if it has the following property:
for a non-isotrivial trigonal curve~$\CC$, one has
$\BM_\CC(\MG)\prec G$
if and only if $\CC$ is Nagata equivalent to
a curve induced from~$\CC_G$.
\enddefinition

\corollary\label{universal.MG}
Any genus zero subgroup $G\subset\MG$ admits a
universal trigonal curve $\CC_G\subset\Sigma_d\to\BB_G\cong\Cp1$
\rom(for some $d$ depending on~$G$\rom)
with type $\tA$ singular fibers only. Such a
universal curve is unique up to
isomorphism.
\endcorollary

\proof
One can take for $\CC_G$ the only trigonal curve determined by the
$j$-invariant $j_G\:\BB_G\to\Cp1$ and the type specification
assigning type~$\bA$ to each singular fiber. In other words,
$\CC_G$ is the maximal trigonal curve determined by the skeleton
$\Sk_G$, \cf. Remark~\ref{rem.group->curve}. Then, the `if' part
in Definition~\ref{def.universal} follows from
Lemmas~\ref{monodromy.induced} and~\ref{monodromy.Nagata}, and the
`only if' part is given by Lemma~\ref{j.factor}.
\endproof

\example\label{ex.ultimate}
The concept of universal curve can be made very explicit in the
case $G=\MG$.
In this case $\tj=j_\CC$.
The `ultimate' universal trigonal curve
$\CC_\MG\subset\Sigma_1$ is the cubic in
the blown-up plane~$\Sigma_1$
given by
$$
\ty^3-3\tx(\tx-1)\ty+2\tx(\tx-1)^2=0,
$$
its $j$-invariant is the identity function $j_\MG(\tx)=\tx$, and its
singular fibers are of types $\tA_0^{**}$ (over $\tx=0$),
$\tA_1^*$ (over $\tx=1$), and $\tA_0^*$ (over $\tx=\infty$).
It is
straightforward that the Weierstra{\ss} equation~\eqref{eq.W} of
any other trigonal curve~$\CC$ is obtained
from the above equation of~$\CC_\MG$ by the substitution
$$
\tx=j_\CC(x),\qquad
\ty=\frac{pq}\Delta y,
$$
see~\eqref{eq.j}, the $y$-substitution corresponding partially to
a sequence of elementary Nagata transformations.
\endexample

\subsection{Subgroups of $\BG3$ and $\tMG$}\label{s.B3.tMG}
Corollary~\ref{universal.MG} has counterparts for the braid
group~$\BG3$ and the extended modular group $\tMG$.

\definition\label{def.universal.B3}
Let $G\subset\BG3$ be a subgroup of genus zero, and let
$m=\depth G$.
A trigonal curve
$\CC_G\subset\Sigma_d\to\BB_G$ is called a \emph{universal
trigonal curve} corresponding to~$G$ if it has the following
property:
for a non-isotrivial trigonal curve~$\CC$ one has
$\BM_\CC(\BG3)\prec G$
if and only if $\CC$ is $m$-Nagata equivalent to
a curve induced from~$\CC_G$.
\enddefinition

\proposition\label{universal.BG3}
Any genus zero subgroup $G\subset\BG3$ admits a universal trigonal
curve\rom; it is unique up
to isomorphism and $m$-Nagata equivalence.
\endproposition

\proof
It suffices to consider the universal curve $\CC_{G'}$
corresponding to the projection $G'$ of~$G$ to~$\MG$, see
Corollary~\ref{universal.MG}, and then change the type
specification to make sure that the new
$\BG3$-valued monodromy
group is a subgroup of~$G$.
\endproof

\Remark\label{rem.proper}
Note that we do \emph{not} assert that the $\BG3$-valued monodromy
group $\BM$ of the curve~$\CC_G$ constructed in the proof
coincides with~$G$. The two groups have the same image in~$\MG$,
but $\depth\BM$ may be a proper multiple of~$m$.
Note, however, that one can
make the two groups coincide by applying an $m$-fold Nagata
transformation at a nonsingular fiber of~$\CC_G$,
introducing an extra
singular fiber with the monodromy $(\Gs_2\Gs_1)^{3m}$.
\endRemark

Given a subgroup $\tG\subset\tMG$, one can consider its
pull-back $G\subset\BG3$ and thus speak about a universal trigonal
curve $\CC_{\tG}$ corresponding to~$\tG$. Note that
$\depth\tG$
equals~$1$ or~$2$ depending on
whether $\tG$ does or, respectively, does not contain $-\id$.

\corollary\label{universal.simple}
Any genus zero subgroup $\tG\subset\tMG$ admits a \emph{simple}
universal trigonal curve $\CC_{\tG}$. If
$-\id\notin\tG$, such a simple curve is unique up to isomorphism.
If, in addition, all singular fibers of $\CC_{\tG}$ are
of Kodaira type~$\I$,
then any \emph{simple} non-isotrivial trigonal curve $\CC$
with $\Im_\CC(\tMG)\prec\tG$ is induced from~$\CC_{\tG}$.
\endcorollary

\proof
It suffices to observe that each series of $2$-Nagata equivalent
singular fibers contains a unique simple one,
see~\eqref{eq.types},
and that any curve
induced from a curve with
Kodaira type~$\I$
singular fibers only is simple, see Lemma~\ref{type.I}.
\endproof

\Remark\label{rem.simple.proper}
The monodromy group of the \emph{simple} universal curve
given by the lemma
may be a proper
subgroup of the pull-back of~$\tG$ in~$\BG3$, \cf.
Subsection~\ref{s.Gamma5} below.
\endRemark

\subsection{Quotients of the fundamental group}\label{s.IK}
We start with a simple lemma.

\lemma\label{I.K}
Let $G\times A\to A$, $(g,a)\mapsto g(a)$, be a group action of a
group~$G$ on a group~$A$, and let $K\subset A$ be a subgroup. Then
the set
$$
\IK_K:=\bigl\{g\in G\bigm|
 \text{\rm$g(a)a\1\in K$ for any $a\in A$}\bigr\}
$$
is a subgroup of~$G$. Furthermore, for each $g\in\IK_K$ one has
$g(K)=K$.
\endlemma

\proof
First, we show that $g(K)\subset K$ for any $g\in\IK_K$. Indeed, if
$k\in K$, then
$$
g(k)=(g(k)k\1)\cdot k\in K
$$
as well. Now, if $g,h\in\IK_K$ and $a\in A$, then
$$
gh(a)a\1=g(h(a)a\1)\cdot g(a)a\1\in K,
$$
\ie, $gh\in\IK_K$, and
$$
g\1(a)a\1=(g(b)b\1)\1\in K,\quad
\text{where}\quad b=g\1(a),
$$
\ie, $g\1\in\IK_K$.
In particular, $g\1(K)\subset K$, hence $g(K)=K$.
\endproof

\definition\label{def.factor}
The affine (projective)
group $\pi:=\piaff\CC$
(respectively, $\piproj\CC$)
of
a trigonal curve $\CC$ is said to \emph{factor} an
epimorphism $\kk\:\AA\onto G$ (or $\kk$ is said to \emph{factor}
through~$\pi$)
if $\kk$ factors through an appropriate
geometric presentation of~$\pi$, \ie, one has
$\kk\:\AA\onto\pi\onto G$ for some geometric presentation
$\AA\onto\pi$ and some epimorphism $\pi\onto G$.
\enddefinition

Fix an epimorphism $\kk\:\AA\onto G$ and let
$K=\Ker\kk\subset\AA$. In view of Lemma~\ref{I.K}, this subgroup
gives rise to a subgroup $\IK_\kk:=\IK_K\subset\BG3$. Due to
Corollary~\ref{groups} (see also
Remark~\ref{rem.geometric}), $\kk$ factors through $\piaff\CC$ if
and only if $\Im_\CC(\BG3)\prec\IK_\kk$. Thus, the following two
statements are immediate consequences of Corollary~\ref{genus.0}
and Proposition~\ref{universal.BG3}, respectively.

\proposition\label{G.genus.0}
An epimorphism $\kk\:\AA\onto G$ factors through the affine group
$\piaff\CC$ of some non-isotrivial
trigonal curve~$\CC$
if and only if $\IK_\kk\subset\BG3$ is a subgroup of
genus zero.
\qed
\endproposition

\proposition\label{G.universal}
Assume that the subgroup $\IK_\kk\subset\BG3$ corresponding to an
epimorphism $\kk\:\AA\onto G$ is of genus zero, and let
$m=\depth\IK_\kk$. Then there exists a unique,
up to isomorphism and $m$-Nagata equivalence,
trigonal curve $\CC_\kk$
with the following property\rom:
the group $\piaff\CC$ of a non-isotrivial trigonal curve $\CC$
factors~$\kk$ if and only if $\CC$ is $m$-Nagata
equivalent to a curve induced from~$\CC_\kk$.
\qed
\endproposition

Any curve $\CC_\kk$ given by Proposition~\ref{G.universal} is called
a \emph{universal} trigonal curve corresponding to~$\kk$. Using
Corollaries~\ref{cor.Nagata} and~\ref{cor.induced},
one arrives at the following statement.

\corollary\label{larger.group}
The group $\piaff\CC$ of a non-isotrivial trigonal curve~$\CC$
factors an epimorphism $\kk\:\AA\onto G$ if
and only if it factors the epimorphism
$$
\AA\onto\piaff{\CC_\kk}/
 \<\text{$[\Gr^m,\Ga]=1$ for any $\Ga\in\AA$}\>
$$
induced from a geometric presentation $\AA\onto\piaff{\CC_\kk}$,
where $m=\depth\IK_\kk$ and
$\CC_\kk$ is \rom(any\rom) universal curve corresponding
to~$\kk$.
\qed
\endcorollary

Due to Corollary~\ref{cor.Nagata}, the quotient
$\piaff{\CC_\kk}\!/[\piaff{\CC_\kk},\<\Gr^m\>]$
in Corollary~\ref{larger.group} does not depend on the choice of a
universal curve.

Finally, assume that $m=\depth\IK_\kk=2$, \ie, $\IK_\kk$ is the
pull-back of a certain subgroup of~$\tMG$ not containing $-\id$.
Then Corollary~\ref{universal.simple} asserts that there is a
unique, up to isomorphism, \emph{simple} universal trigonal
curve~$\CC_\kk$. Furthermore, combining this observation with
Corollary~\ref{cor.induced},
one has the following statement.

\corollary\label{simple.group}
Assume that $\depth\IK_\kk=2$ and that the simple universal
curve~$\CC_\kk$ has Kodaira type~$\I$ singular fibers only. Then
the group $\piaff\CC$ of a \emph{simple}
non-isotrivial
trigonal curve~$\CC$
factors~$\kk$ if and only if it factors a geometric presentation
$\AA\onto\piaff{\CC_\kk}$ of the group of~$\CC_\kk$.
\qed
\endcorollary

\section{Uniform dihedral quotients\label{S.DG}}

The principal result of this section is the proof of
Theorems~\ref{th.D} and~\ref{th.reducible}. In
Subsection~\ref{s.isotrivial}, we treat in details the case of
isotrivial curves, which form a very spacial subclass, mostly
excluded from the consideration in the previous sections.

\subsection{Quotients to be considered}\label{s.quotients}
Consider the epimorphism $\deg_2\:\AA\onto\CG2$ given by
$\Ga\mapsto\deg\Ga\bmod2$.
By the Reidemeister--Schreier
algorithm (see, \eg,~\cite{MKS}), the kernel $K:=\Ker\deg_2$ is
freely generated by the elements $u:=\Ga_1^2$, $v_1=\Ga_1\Ga_2$,
$v_2:=\Ga_2\Ga_1$, $w_1:=\Ga_1\Ga_3$, and $w_2:=\Ga_3\Ga_1$, and
the conjugation~$t$ by the generator of $\CG2=\AA/K$ acts on the
abelianization $K/[K,K]$ \via\
$u\leftrightarrow u$, $v_1\leftrightarrow v_2$,
$w_1\leftrightarrow w_2$.

The maximal generalized dihedral group through which $\deg_2$
factors is
$\GDG(\CH)$, where $\CH$ is the quotient of
$K/[K,K]$ by the subgroup $\Im(t+\id)$ and the square of (any)
representative of the generator of~$\CG2$. (Modulo $\Im(t+\id)$,
this square does not depend on the choice of a representative.)
The computation above shows that $\CH=\Z a\oplus\Z b$, where
$a:=v_2=-v_1$ and $b:=w_1=-w_2$. In other words, $\CH$ is the
abelianization of the group
$\Ker{\deg_2}/\<\Ga_1^2,\Ga_2^2,\Ga_3^2\>$.

It follows that any generalized dihedral quotient of~$\AA$ that
factors~$\deg_2$ is of the form
$$
\AA\onto\GDG(\CH)\onto\GDG(\CH/H)\eqtag\label{eq.uniform}
$$
for some subgroup $H\subset\CH$.

\definition\label{def.uniform}
An epimorphism $\kk\:\AA\onto\GDG(\CH/H)$, $H\subset\CH$,
as in~\eqref{eq.uniform} is
called a \emph{uniform}
dihedral quotient of~$\AA$.
If $\kk$ factors through
a geometric presentation $\AA\onto\piaff\CC$ of the
fundamental
group
of a
trigonal curve~$\CC$, then $\GDG(\CH/H)$ is also said to be a
\emph{uniform}
dihedral quotient of $\piaff\CC$, and
$\CC$ itself is said to admit a \emph{uniform}
$\GDG(\CH/H)$-covering.
\enddefinition

We identify the group~$\CH=\Z a\oplus\Z b$
obtained in the previous paragraph
with the group~$\CH$
introduced in Subsection~\ref{s.groups}.
The kernel of the epimorphism $\AA\onto\GDG(\CH)$ is
$\BG3$-invariant; hence, the $\BG3$-action on~$\AA$ induces a
certain action on $\GDG(\CH)$.

\lemma\label{induced.action}
The induced $\BG3$-action on $\GDG(\CH)=\CH\rtimes\CG2$
splits into the product
of the canonical representation of~$\tMG$ on~$\CH$ and the trivial
action on~$\CG2$.
\endlemma

\proof
The induced action on~$\CH$ is well known; it can
easily
be found
by a direct computation, \cf. Remark~\ref{modular.bm}.
The splitting follows from the fact that
the generator of~$\CG2$ admits an invariant representative,
namely~$\Gr$.
\endproof

\corollary\label{cor.CQ}
Any uniform dihedral quotient of
the group~$\piaff\CC$ factors through the
\emph{maximal uniform quotient} $\piaff\CC\onto\GDG(\CQ)$, where
$\CQ=\CH/H$ and $H$ is the sum of the images $\Im(g-\id)$ over all
$g\in\BM_\CC(\tMG)$.
\qed
\endcorollary

\Remark\label{rem.BG2}
Similarly, for the standard $\BG{n}$-action on the free
group~$F_n$, one can consider the maximal invariant dihedral
quotient $F_n\onto\GDG(\CH)$, $\CH\cong\bigoplus_{n-1}\Z$.
Lemma~\ref{induced.action} is in a sharp contrast with the case of
$n$ even. For example, if $n=2$, the induced action on
$\CH\cong\Z$ is trivial, whereas the action on $\GDG(\CH)$ is not.
For this reason, the maximal uniform dihedral quotient of the
fundamental group of a hyperelliptic (`bigonal') curve
in~$\Sigma_d$ is not controlled by the always trivial action of
the monodromy group of the curve on the kernel~$\CH$, \cf.
Subsection~\ref{s.groups} below. If the monodromy group is
generated by~$\Gs_1^s$, $s>0$, the maximal dihedral quotient is
easily found to be $\DG{2s}$. In particular, for each integer
$s>0$, there exists a reducible trigonal curve in $\Sigma_{2s}$
whose fundamental group has a (non-uniform) quotient
$\DG{2s}$.
\endRemark

Next statement justifies
our interest in uniform dihedral quotients;
for more motivation, see Section~\ref{S.appl} below.

\proposition\label{all.uniform}
Any generalized dihedral quotient of the fundamental group
$\piaff\CC$ of an irreducible trigonal curve~$\CC$ is uniform.
\endproposition

\proof
If $\CC$ is irreducible, the abelianization of $\piaff\CC$
equals~$\Z$; hence $\piaff\CC$ admits a unique epimorphism
to~$\CG2$, and this epimorphism factors~$\deg_2$. On the
other hand, $\GDG(\CH)$ is the maximal generalized dihedral group
through which $\deg_2$ factors.
\endproof

\subsection{The monodromy groups}\label{s.subgroups}
Consider a uniform
generalized
dihedral quotient $\kk\:\AA\onto\GDG(\CH/H)$, $H\subset\CH$.
Due to Lemma~\ref{induced.action},
the subgroup $\IK_\kk\subset\BG3$ introduced in
Subsection~\ref{s.IK} is the pullback of the subgroup
$$
\IK_H:=\bigl\{g\in\tMG\bigm|\Im(g-\id)\subset H\bigr\}\subset\tMG.
$$
Up to the action of~$\tMG$, one has
$H=\Z(ma)\oplus\Z(nb)\subset\CH=\Z a\oplus\Z b$, where $m$, $n$
are nonnegative integers and
either
$m=n=0$ or $m\ne0$ and $m\divides|n$.
We consider separately the following three cases.

\subsubsection{The case $m=n=0$}\label{ss.m=n=0}
In this case $H=0$ and $\IK_H=\{1\}$. Hence, $\IK_\kk\subset\BG3$
is the central cyclic subgroup generated by $(\Gs_2\Gs_1)^6$.

\subsubsection{The case $m\ne0$, $n=0$}\label{ss.n=0}
In this case, $\IK_H$ is the infinite cyclic subgroup generated by
(the image of) $\Gs_1^m$.
Note that both $[\tMG:\IK_H]$ and $[\BG3:\IK_\kk]$ are infinite.

\subsubsection{The case $m,n\ne0$, $m\divides|n$}\label{ss.m|n}
In this case, $\IK_H$ is the subgroup
$$
\tMG_m(n):=\left\{\bmatrix a&b\\c&d\endbmatrix\in\tMG\left|
 \bmatrix a&b\\c&d\endbmatrix=\bmatrix 1&*\\0&1\endbmatrix\bmod n,\
 b=0\bmod m
\right\}\rlap.\right.
$$
(We use the notation of~\cite{Cox}.)
The image of $\tMG_m(n)$ in~$\MG$ is
$\MG_m(n):=\MG_1(n)\cap\MG(m)$; note that $\MG_1(n)$ is consistent
with the conventional notation. Note also that $\MG_n(n)$ is
the principal congruence subgroup $\MG(n)$ and
$\MG_m(n)\supset\MG(n)$ for any $m\divides|n$.
In particular, all
groups $\IK_H$ obtained in this case are congruence subgroups.

In the sequel, we use repeatedly the following obvious
observation: if $g\in\tMG$ and $\det(g-\id)=d\ne0$, then
$\Im(g-\id)\subset\CH$ is a subgroup of index~$d$.

\lemma\label{MG.torsion.free}
Unless
$m=1$ and $n\le3$,
the group $\MG_m(n)$ is
torsion
free.
\endlemma

\proof
Any torsion element of~$\MG$ is conjugate to either $\X^{\pm1}$
or~$\Y$, and for any lift $g\in\tMG$ of such an element the
determinant $\det(g-\id)$ takes value~$1$, $2$, or~$3$.
\endproof

\lemma\label{MG.-id}
The group $\tMG_m(n)$ contains $-\id$ if and only if $n\le2$.
\endlemma

\proof
One has $\Im[-\id-\id]=2\CH$.
\endproof

Recall that \emph{parabolic} elements of~$\tMG$ are those of the
form $\pm g$, where $g\in\tMG$ is a unipotent matrix, \ie,
$(g-\id)^2=0$. Any unipotent element is conjugate to a power of
$\X\Y$ (the image of~$\Gs_1$). Up to~$\pm\id$, any torsion free
genus zero subgroup is generated by parabolic elements.

\lemma\label{MG.parabolic}
If $n>2$ and $(m,n)\ne(1,4)$, then any parabolic element
of~$\tMG_m(n)$ is unipotent.
\endlemma

\proof
If $g\in\tMG$ is unipotent, then $\det(-g-\id)=4$.
\endproof

\lemma\label{MG(1,4)}
The group $\tMG_1(4)$ has three conjugacy classes of indivisible
parabolic elements, those of the images of~$\Gs_1$, $\Gs_2^4$, and
$\Gs_1\Gs_2^4$. The first two are unipotent, the last one is not.
\endlemma

\proof
The proof is a direct computation, see, \eg,
Figure~\ref{fig.curves.r}(c) below, where the skeleton
$\MG\!/\MG_1(4)$ is shown.
\endproof

\subsection{Proof of Theorem~\ref{th.reducible}}\label{proof.reducible}
Consider the
epimorphism
$\BG3\onto\SG3=\smash{\tMG/\tMG(2)}$, and
let $\BM_\CC(\SG3)\subset\SG3$ be the image of $\BM_\CC(\BG3)$.
Geometrically, $\BM_\CC(\SG3)$ is the monodromy group of the
ramified covering $\CC\to\Cp1$. Hence, $\CC$ is reducible if and
only if $\BM_\CC(\SG3)$ is not transitive, \ie,
$\BM_\CC(\tMG)\prec\tMG_1(2)$, and $C$ splits into three components
if and only if $\BM_\CC(\SG3)=\{1\}$, \ie,
$\BM_\CC(\tMG)\subset\tMG(2)$. By~\ref{ss.m|n} above, these
conditions on $\BM_\CC(\tMG)$
are equivalent to the existence of uniform
$\GDG(\CG2)$- and $\GDG(\CG2\oplus\CG2)$-coverings, respectively.

According to Proposition~\ref{all.uniform}, any dihedral covering of
an irreducible trigonal curve is uniform. Hence, any curve
admitting any $\GDG(\CG2)$-covering is reducible.

The group $\GDG(\CG2\oplus\CG2)\cong\CG2\oplus\CG2\oplus\CG2$ is
abelian with three generators. Hence, if $\piaff\CC$ factors to
$\GDG(\CG2\oplus\CG2)$, its abelianization has at least three
generators; from the Poincar\'e--Lefschetz duality (applied to the
union $\CC\cup\EE\cup F\subset\Sigma_d$) it follows that $\CC$ must
have at least three components.
\qed

\Remark\label{rem.S3.torus}
In the proof of
Lemma~\ref{torus.monodromy}, we considered a
trigonal curve $\CC_x\subset\Sigma_1$ with the set of singular
fibers $\tA_1^*\oplus3\tA_0^*$. Any such curve can be obtained
by a perturbation from a curve $\CC_x'\subset\Sigma_1$ with
the singular fibers $\tA_1^*\oplus\tA_0^{**}\oplus\tA_0^*$ (a
nodal plane cubic projected from the point of intersection of the
tangents to one of the branches at the node and one of the
inflection points). The skeleton of $\CC_x'$ is the graph
$\mathord\circ{\joinrel\relbar\joinrel\relbar\joinrel}\mathord\bullet$
corresponding to the full modular group~$\MG$. Hence, the
$\MG$-valued monodromy groups of both curves are~$\MG$ and,
combining with the epimorphism $\MG\onto\SG3$ above, one concludes
that the $\SG3$-valued groups are~$\SG3$, as stated.
\endRemark

\subsection{Isotrivial curves}\label{s.isotrivial}
Isotrivial trigonal curves are easily classified and their
monodromy groups are easily computed.
Depending on the constant $j_\CC$, one can
distinguish the following three cases.

\subsubsection{The case $j_\CC\equiv0$}\label{ss.j=0}
One has $p(x)\equiv0$ in~\eqref{eq.W} and~\eqref{eq.j} and
the Weierstra{\ss} equation takes the form
$y^3+2q(x)=0$. The monodromy group $\BM_\CC(\BG3)$ is the cyclic
group generated by $(\Gs_2\Gs_1)^r$, where $r$ is the greatest
common divisor of the multiplicities of the roots of $q(x)$.
Hence, $\BM_\CC(\tMG)$ is generated by $(-\X)^r$, and the maximal
uniform dihedral quotient that $\piaff\CC$ may have is
$\GDG(\CH/\Im[(-\X)^r-\id])$. Note that
$\det[(-\X)^{\pm1}-\id]=1$ and $\det[(-\X)^{\pm2}-\id]=3$,
whereas $(-\X)^3=-\id$.
Summarizing, one has the following statements:
\roster
\item\local{j=0,r=0}
if $r=0\bmod6$, the maximal
uniform
quotient is
$\piaff\CC\onto\GDG(\CH)$;
\item\local{j=0,r=3}
if $r=3\bmod6$, the maximal
uniform
quotient is
$\piaff\CC\onto\GDG(\CG2\oplus\CG2)$;
\item\local{j=0,r=1}
if $r=\pm1\bmod6$, then $\piaff\CC$ admits no nontrivial dihedral
quotients;
\item\local{j=0,r=2}
if $r=\pm2\bmod6$, the maximal generalized
dihedral quotient of $\piaff\CC$ is $\DG6$.
\endroster
In a sense, the curves as in~\loccit{j=0,r=0} and~\loccit{j=0,r=3}
are degenerate and can be regarded as a special case
of~\ref{ss.j=const} below.
In cases~\loccit{j=0,r=1} and~\loccit{j=0,r=2}, the curves are
irreducible and all their generalized dihedral coverings are uniform.

Observe that, in case~\loccit{j=0,r=2}, adding an extra
relation $(\Gs_2\Gs_1)^6=\id$, \ie, making~$\Gr^2$ a central
element,
one obtains an epimorphism
$\piaff\CC\onto\BG3$ (given by $\Ga_1,\Ga_3\mapsto\Gs_1$,
$\Ga_2\mapsto\Gs_2$).
Observe also that, in this case, all roots of~$q$ have even
multiplicities; hence, one has $q=\bar q^2$ for some polynomial
$\bar q(x)$ and the equation of the curve has the form
$y^3+(\sqrt2\bar q)^2=0$, \ie, the curve is of torus
type, see Subsection~\ref{s.torus}.

\subsubsection{The case $j_\CC\equiv1$}\label{ss.j=1}
One has $q(x)\equiv0$ in~\eqref{eq.W} and~\eqref{eq.j} and
the Weierstra{\ss} equation takes the form
$y(y^2+3p(x))=0$. The monodromy group $\BM_\CC(\BG3)$ is the cyclic
group generated by $(\Gs_2\Gs_1^2)^r$, where $r$ is the greatest
common divisor of the multiplicities of the roots of $p(x)$.
Hence, $\BM_\CC(\tMG)$ is generated by $(-\Y)^r$, and the maximal
uniform dihedral quotient that $\piaff\CC$ may have is
$\GDG(\CH/\Im[(-\Y)^r-\id])$. Note that
$\det[(-\Y)^{\pm1}-\id]=2$,
whereas $(-\Y)^2=-\id$.
Thus, one has:
\roster
\item\local{j=1,r=0}
if $r=0\bmod4$, the maximal
uniform
quotient is
$\piaff\CC\onto\GDG(\CH)$;
\item\local{j=1,r=2}
if $r=2\bmod4$, the maximal
uniform
quotient is
$\piaff\CC\onto\GDG(\CG2\oplus\CG2)$;
\item\local{j=1,r=1}
if $r=1\bmod2$, the maximal
uniform
quotient is
$\piaff\CC\onto\DG4$.
\endroster
As in~\ref{ss.j=0} above,
the curves as in~\loccit{j=1,r=0} and~\loccit{j=1,r=2} splitting
into three components
can be regarded as a special case of~\ref{ss.j=const} below.

\subsubsection{The case $j_\CC=\const\ne0,1$}\label{ss.j=const}
In this case, one has $p^3/q^2=\const\ne0$. Hence, there is a
polynomial $s(x)$ such that $p=\Ga s^2$ and $q=\Gb s^3$,
$\Ga,\Gb=\const$.
It is straightforward that the curve is Nagata equivalent to the
union of three generatrices of the form $y=\const$ in $\Sigma_0$.
(Occasionally, a curve~$\CC$ of this form may have $j_\CC=0$
or~$1$, corresponding to cases~\iref{ss.j=0}{j=0,r=0},
\ditto{j=0,r=3} and~\iref{ss.j=1}{j=1,r=0},
\ditto{j=1,r=2}, respectively,
see remarks
in the corresponding sections.)
Each $m$-fold
Nagata transformation results in a type~$\tJ_{m,0}$ singular
fiber. Hence, the set of singular fibers of~$\CC$ is of the form
$\bigoplus\tJ_{m_i,0}$, and
the monodromy group
$\BM_\CC(\BG3)$ is
the cyclic group generated by $(\Gs_2\Gs_1)^{3r}$, where
$r=\gcd(m_i)$. Thus, the maximal uniform dihedral quotient of
$\piaff\CC$ is $\GDG(\CH/\Im[(-\id)^r-\id])$. One has:
\roster
\item\local{j=const,r=0}
if $r=0\bmod2$, the maximal
uniform
quotient is
$\piaff\CC\onto\GDG(\CH)$;
\item\local{j=const,r=1}
if $r=1\bmod2$, the maximal
uniform
quotient is
$\piaff\CC\onto\GDG(\CG2\oplus\CG2)$.
\endroster

\definition\label{def.trivial}
A trigonal curve as in~\iref{ss.j=0}{j=0,r=0},
\iref{ss.j=1}{j=1,r=0}, or~\iref{ss.j=const}{j=const,r=0} above is
called \emph{trivial}.
\enddefinition

Trivial
curves are those
$2$-Nagata equivalent to the union of three `horizontal'
generatrices in $\Sigma_0=\Cp1\times\Cp1$. All such curves
split into three components, their $\tMG$-valued monodromy groups
are trivial,
and their fundamental groups admit infinite uniform dihedral
quotients $\piaff\CC\onto\GDG(\Z\oplus\Z)$. A trivial curve
$\CC\subset\Sigma_d$ can also be characterized as follows:
$d$ is even and
the
relatively minimal model of
the double covering of~$\Sigma_d$
ramified at $\CC\cup\EE$ is a trivial elliptic surface
$\text{(elliptic curve)}\times\Cp1$.
The latter can also be described as the only irregular elliptic
surfaces (over a rational base) or the only elliptic surfaces
without singular fibers.

\Remark\label{rem.isotrivial}
It follows that the
conclusion of Theorem~\ref{th.D} holds for any nontrivial
isotrivial curve.
\endRemark

\subsection{Proof of Theorem~\ref{th.D} and
Corollary~\ref{cor.irreducible}}\label{proof.D}
The case of isotrivial curves is considered in
Subsection~\ref{s.isotrivial}, see Remark~\ref{rem.isotrivial}.

Assume that $\CC$ is non-isotrivial.
Due to Proposition~\ref{G.genus.0},
a uniform dihedral quotient $\AA\onto\GDG(\CH/H)$
factors through
$\piaff\CC$
(for some non-isotrivial curve~$\CC$)
if and only if
the subgroup $\IK_H\subset\tMG$ introduced in
Subsection~\ref{s.subgroups}
is of genus zero.
Since the cases considered in~\ref{ss.m=n=0} and~\ref{ss.n=0}
give rise to subgroups of infinite index,
$\IK_H$ must be of the form $\tMG_m(n)$ for some
$m,n>0$, $m\divides|n$, the resulting quotient being $\DG{2n}$ if
$m=1$ or $\GDG(\CG{m}\oplus\CG{n})$ if $m>1$.
The genera of $\MG_m(n)$ are computed in~\cite{Cox}, and the
subgroups of genus zero are precisely those given by
Theorem~\ref{th.D}. Alternatively, one can observe that, with the
exception of $\MG_1(1)=\MG$, $\MG_1(2)$, and $\MG_1(3)$,
all groups $\MG_m(n)$ are
torsion free, see Lemma~\ref{MG.torsion.free},
and refer to the list of torsion free genus zero congruence
subgroups found in~\cite{Sebbar}.

Corollary~\ref{cor.irreducible} follows
from Theorems~\ref{th.D},~\ref{th.reducible} and
Proposition~\ref{all.uniform}.
\qed

\section{Properties of universal curves\label{S.groups}}

We construct the minimal universal curves corresponding to
uniform dihedral coverings,
discuss their geometric properties, and, for the irreducible
curves, compute their
fundamental groups.
The results are applied to illustrate
Statements \iref{spec}{spec.geometry} and~\ditto{spec.larger} in
the introduction
(see, \eg,
Corollaries~\ref{cor.Phi6} and~\ref{cor.Phi10})
and, in
particular, to prove Theorem~\ref{th.Oka}.
In Remark~\ref{rem.Delta2},
we announce a few
further results concerning the Alexander polynomial.

\subsection{The universal curves corresponding to dihedral
quotients}
Consider a uniform generalized dihedral quotient
$\kk\:\AA\to\GDG(\CH/H)$ and let $\IK_H\subset\tMG$ be the
subgroup introduced in Subsection~\ref{s.subgroups}. Due to
Corollary~\ref{universal.simple}, there is a simple universal
curve~$\CC_H$ corresponding to~$\IK_H$. If $\IK_H\not\ni-\id$,
this curve is unique up to isomorphism; otherwise, $\CC_H$ is
defined up to Nagata equivalence, and one can choose all singular
fibers of~$\CC_H$ to be of type~$\tA$.

All universal curves are maximal. Their skeletons
$\Sk_H\subset S^2$, marked by the
isomorphism type of the quotient $\CH/H$, are shown in
Figures~\ref{fig.curves} (irreducible curves)
and~\ref{fig.curves.r} (reducible curves). We omit the trivial
case $H=\CH$; the corresponding ultimate trigonal curve is
discussed in Example~\ref{ex.ultimate}. In the figures, the
sphere~$S^2$ is represented as the quotient space $D/\partial D$,
where $D$ is the disk bounded by the grey dotted circle.

To complete the description of the universal curves, we need to
determine the type
specifications.
This is done in~\ref{ss.2}--\ref{ss.others} below.

\midinsert
\def\tl#1#2#3{\relax\hbox to0pt{\hss(#1)\enspace
$#3$\hss}}
\centerline{\vbox{\halign{\hss#\hss&&\qquad\hss#\hss\cr
 \cpic{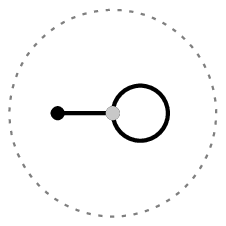}&\cpic{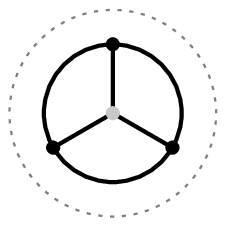}&\cpic{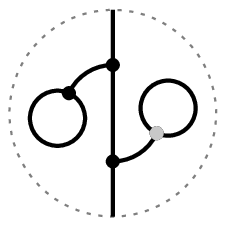}&\cpic{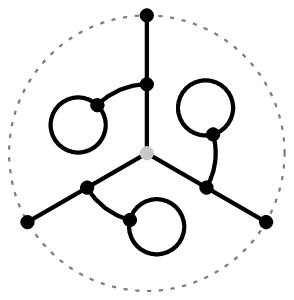}\cr
 \noalign{\medskip}
 \tl a{\MGd3}{\CG3}&\tl b{\MG(3)}{\CG3\oplus\CG3}&
 \tl c{\MGd5}{\CG5}&\tl d{\MGd7}{\CG7}\cr}}}
\bigskip
\centerline{\vbox{\halign{\hss#\hss&&\qquad\qquad\hss#\hss\cr
 \cpic{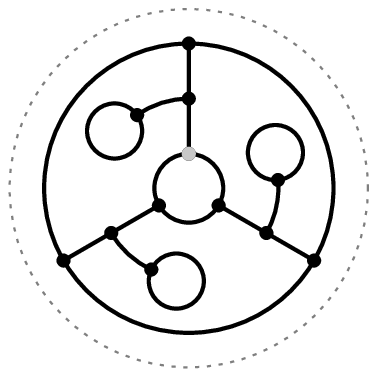}&\cpic{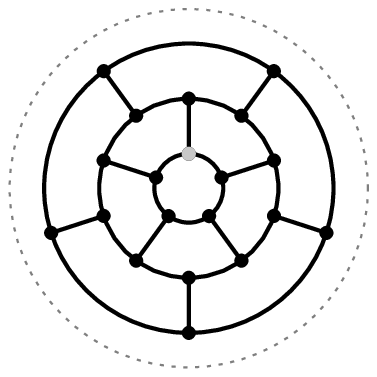}\cr
 \noalign{\medskip}
 \tl e{\MGd9}{\CG9}&\tl f{\MG(5)}{\CG5\oplus\CG5}\cr}}}
\figure\label{fig.curves}
Irreducible universal curves admitting dihedral coverings
\endfigure
\endinsert

\midinsert
\def\tl#1#2#3{\relax\hbox to0pt{\hss(#1)\enspace
$#3$\hss}}
\centerline{\vbox{\halign{\hss#\hss&&\quad\ \ \hss#\hss\cr
 \cpic{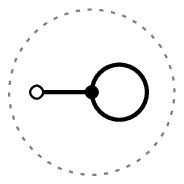}&\cpic{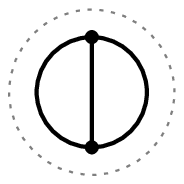}&\cpic{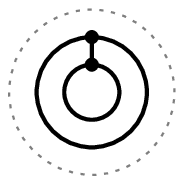}&\cpic{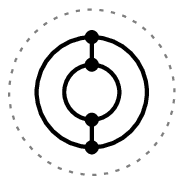}&\cpic{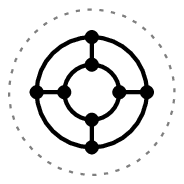}\cr
 \noalign{\medskip}
 \tl a{\MGd2}{\CG2}&\tl b{\MG(2)}{\CG2\oplus\CG2}&\tl c{\MGd4}{\CG4}&
 \tl d{\MG_2(4)}{\CG2\oplus\CG4}&\tl e{\MG(4)}{\CG4\oplus\CG4}\cr}}}
\bigskip
\centerline{\vbox{\halign{\hss#\hss&&\qquad\hss#\hss\cr
 \cpic{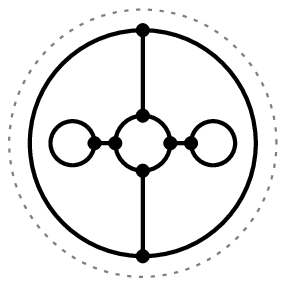}&\cpic{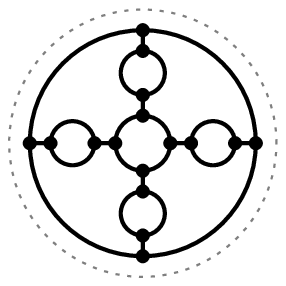}&\cpic{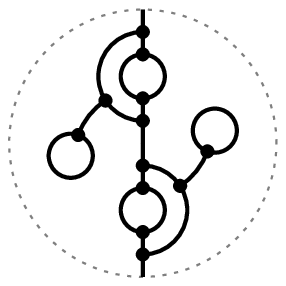}\cr
 \noalign{\medskip}
 \tl f{\MGd8}{\CG8}&\tl g{\MG_2(8)}{\CG2\oplus\CG8}&\tl h{\MGd{10}}{\CG{10}}\cr}}}
\bigskip
\centerline{\vbox{\halign{\hss#\hss&&\qquad\hss#\hss\cr
 \cpic{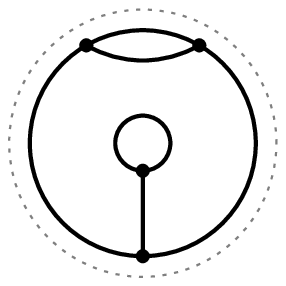}&\cpic{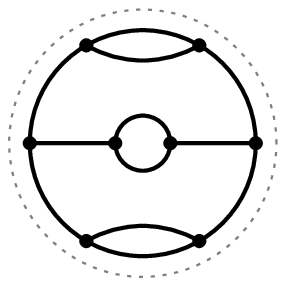}&\cpic{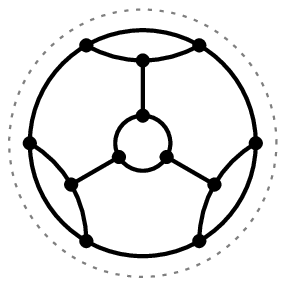}\cr
 \noalign{\medskip}
 \tl i{\MGd6}{\CG6}&\tl j{\MG_2(6)}{\CG2\oplus\CG6}&\tl k{\MG_3(6)}{\CG3\oplus\CG6}\cr}}}
\figure\label{fig.curves.r}
Reducible universal curves admitting dihedral coverings
\endfigure
\endinsert

\subsubsection{The groups~$\tMG_1(2)$ and~$\tMG(2)$}\label{ss.2}
The skeletons are shown in Figures~\ref{fig.curves.r}(a) and~(b),
respectively. Both groups contain $-\id$, see Lemma~\ref{MG.-id},
hence all singular
fibers can be chosen of type~$\tA$. Under this convention, both
curves are cubics in the blown up plane~$\Sigma_1$. The sets of
singular fibers are $\tA_1^*\oplus\tA_1\oplus\tA_0^*$ for
$\tMG_1(2)$ and $3\tA_1$ for $\tMG(2)$. In view of
Theorem~\ref{th.reducible}, the two curves can be regarded as,
respectively,
the universal reducible trigonal curve and
the universal trigonal curve split into three components.

\subsubsection{The group~$\tMG_1(4)$,
Figure~\ref{fig.curves.r}(c)}\label{ss.4}
The group does not contain~$-\id$, hence the curve is unique.
Due to Lemma~\ref{MG(1,4)}, the type specification is
$\tA_3\oplus\tD_5\oplus\tA_0^*$.

\subsubsection{The group~$\tMG_1(3)$,
Figure~\ref{fig.curves}(a)}\label{ss.3}
The group contains two conjugacy classes of torsion elements (the
monovalent vertex in the figure),
which must be $\X^{\pm1}$, as $\det(-\X^{\pm1}-\id)=1$.
All parabolic elements are unipotent, see
Lemma~\ref{MG.parabolic}, corresponding to Kodaira type~$\I$
singular fibers.
Hence, the curve belongs to~$\Sigma_2$ and its type specification
is $\tE_6\oplus\tA_2\oplus\tA_0^*$. This curve can be given by the
affine equation
$y^3+(y+x)^2=0$;
hence it is of torus type.

\subsubsection{Other groups}\label{ss.others}
All other groups are torsion free, see
Lemma~\ref{MG.torsion.free}, do not contain~$-\id$, see
Lemma~\ref{MG.-id}, and all their parabolic elements are
unipotent, see Lemma~\ref{MG.parabolic}. Hence, all their singular
fibers are of Kodaira type~$\I$. The curves belong to~$\Sigma_d$,
where $2d$ is the number of vertices in the corresponding
skeleton.

\subsection{The Alexander polynomial}\label{s.Alexander}
In this subsection, we remind the definition and basic properties
of the Alexander polynomial of an algebraic curve. For more
details, see A.~Libgober~\cite{Libgober}.

Consider a group~$G$ and a distinguished epimorphism
$\kk\:G\onto\Z$. Denote by~$\AM_\kk$ the abelianization of the
kernel $\Ker\kk$. The conjugation by the generator of~$\Z$ induces
an automorphism $t\:\AM_\kk\to\AM_\kk$,
thus making $\AM_\kk$ a module over the
ring $\Lambda:=\Z[t^{\pm1}]$ of Laurent polynomials in~$t$.
The abelian
group~$\AM_\kk$, regarded as a $\Lambda$-module, is called the
\emph{Alexander module} of~$G$ (more precisely, of~$\kk$).

If $G$ is finitely generated, the Reidemeister--Schreier algorithm
(see, \eg,~\cite{MKS}) asserts that $\AM_\kk$ is finitely generated
over~$\Lambda$. Hence, $\AM_\kk\otimes\C$ is finitely generated over
the principal ideal domain $\Lambda_\C:=\C[t^{\pm1}]$. If
$\AM_\kk\otimes\C$ happens to be a torsion module, its order
$\Delta_\kk(t)\in\Lambda_\C$ is called the \emph{Alexander
polynomial} of~$G$. It is defined up to units, \ie, up to
multiplication by $ct^n$, $c\in\C^*$, $n\in\Z$. Alternatively,
under the assumptions $\AM_\kk\otimes\C$ is a finite dimensional
vector space with an operator~$t$, and $\Delta_\kk$ can be defined
as the characteristic polynomial of~$t$. Note that, in fact, $t$
acts on the finitely generated free abelian group $\AM_\kk/\Tors$.
Hence, properly normalized, $\Delta_\kk$ is defined over~$\Z$.

Any epimorphism $\Gr\:G'\onto G$ of groups induces an
epimorphism
$\AM_{\kk\circ\Gr}\onto\AM_\kk$ of their Alexander modules,
which is $\Lambda$-linear. Hence
one has $\Delta_\kk\divides|\Delta_{\kk\circ\Gr}$, provided that
$\Delta_{\kk\circ\Gr}$ is defined.

Given a trigonal curve $\CC\subset\Sigma_d$, take for
$\kk\:\piaff\CC\onto\Z$ the quotient of the degree
epimorphism $\deg\:\AA\onto\ZZ$. The
module $\AM_\CC:=\AM_\kk$ is called the \emph{Alexander module}
of~$\CC$. Note that the central element~$\Gr^d$, see
Corollary~\ref{commutants}, is mapped to $3d\in\Z$; hence, $\AM_\CC$
is annihilated by $t^{3d}-1$. Thus, unless $\CC$ is a trivial
curve in~$\Sigma_0$, the Alexander polynomial
$\Delta_\CC(t):=\Delta_\kk(t)$ is well defined; moreover, it is a
product of cyclotomic polynomials~$\Phi_m$ with $m\divides|3d$.
Furthermore, since $\AM_\AA\cong\Lambda\oplus\Lambda$, each
cyclotomic factor appears in $\Delta_\CC$ at most twice.

The following statement is essentially contained in
O.~Zariski~\cite{Zariski}.

\lemma\label{Delta.p}
If a trigonal curve $\CC$ is irreducible, then $\Delta_\CC(t)$ is
prime to $\Phi_1$ and to $\Phi_m$ for any prime power~$m$.
\qed
\endlemma

\lemma\label{Delta.2p}
Let $m=p^a$ be a prime power. If
$\Phi_{2m}(t)\divides|\Delta_\CC(t)$,
then $\CC$ admits a uniform $\DG{2p}$-covering. In particular, one
has $p=2$, $3$, $5$, or~$7$.
\endlemma

\proof
If $p=2$, the statement follows from Lemma~\ref{Delta.p}, which
asserts that the curve is reducible,
and Theorem~\ref{th.reducible}.

Assume~$p$ odd.
By the assumption, $\AM':=\AM_\CC/\Phi_{2m}$ is an infinite abelian
group; hence one has $\AM'':=\Hom_\Z(\AM',\CG{p})\ne0$. The
$\Z$-action on~$\AM''$ factors to an action of
$\CG{2m}=\CG2\oplus\CG{m}$, the $\CG2$ summand
acting \via~$-\id$. The action of the
$p$-group $\CG{m}$ on the $p$-group~$\AM''$ has an invariant element
$\phi\ne0$,
and the semidirect product
$(\AM'\!/\Ker\phi)\rtimes\CG2\cong\DG{2p}$
is the desired dihedral quotient.
\endproof

\Remark
At present, I do not know whether $\Phi_{14}$ can appear as a
factor of the Alexander polynomial of a trigonal curve. The
factors $\Phi_6$, $\Phi_6^2$, $\Phi_{10}$, and $\Phi_{10}^2$ do
appear, see examples below.
\endRemark

\Remark
One has $\CH=\AM_\AA/(t+1)$ and, in view of~\eqref{eq.uniform},
the maximal uniform dihedral quotient of $\piaff\CC$ is
$\GDG(\CQ)$, where $\CQ=\AM_\CC/(t+1)$.
\endRemark

\subsection{Proof of Theorem~\ref{th.Oka}}\label{proof.Oka}
The implications \itemref{th.Oka}{Oka.B}$\implies$\ditto{Oka.D}
and \ditto{Oka.B}$\implies$\ditto{Oka.Delta} are obvious,
the implication \ditto{Oka.torus}$\implies$\ditto{Oka.D}
is given by Lemma~\ref{torus.monodromy},
and the implication \ditto{Oka.Delta}$\implies$\ditto{Oka.D}
is given by Lemma~\ref{Delta.2p}. It remains to prove that
\ditto{Oka.D}$\implies$\ditto{Oka.B},~\ditto{Oka.torus}.

If $\CC$ is isotrivial, the statement is contained in
Subsection~\ref{s.isotrivial}, the last paragraph
of~\ref{ss.j=0}. If $\CC$ is not isotrivial, it is $2$-Nagata
equivalent to a curve induced from the universal curve~$\CC_3$
described in~\ref{ss.3}. The latter is of torus type,
hence so is~$\CC$, see
Lemma~\ref{torus.invariant}. The braid monodromy of~$\CC_3$ is
found using~\cite{degt.kplets} and the skeleton shown in
Figure~\ref{fig.curves}(a); the monodromy group is generated by
$\Gs_1$ and~$\Gs_2^3$ (the two regions in the figure) and the
monodromy at infinity $(\Gs_2\Gs_1)^6$, see Lemma~\ref{infty}.
The resulting affine group
given by Corollary~\ref{groups} is isomorphic to~$\BG3$,
with~$\Gr^2$ mapped to $(\Gs_2\Gs_1)^6$, and
Corollary~\ref{larger.group} implies that $\piproj\CC$ factors to
$\BG3/(\Gs_2\Gs_1)^6\cong\CG2*\CG3$.
\qed

\Remark\label{rem.Oka}
It follows from the proof that Theorem~\ref{th.Oka} extends to
reducible curves as well, provided that in
Item~\itemref{th.Oka}{Oka.D}
one speaks about uniform quotients to~$\DG6$ and in
Item~\ditto{Oka.B},
about quotients $\piproj\CC\onto\CG2*\CG3$
that further factor to the uniform epimorphism
$\piproj\CC\onto\CG2$.
\endRemark

\subsection{Fundamental groups of universal curves}
In the rest of this section, we compute the fundamental groups of
the other irreducible universal curves. Our principal goal is
exploring more examples illustrating Speculation~\ref{spec}.
Besides, this computation can also serve as a proof of
the fact that the skeletons shown in Figure~\ref{fig.curves} do
indeed represent the universal curves, as for each
admissible pair $(m,n)$, the
corresponding skeleton has the correct number of edges
$[\MG:\MG_m(n)]$ and the resulting curve is irreducible and admits
a $\GDG(\CG{m}\oplus\CG{n})$-covering.

The groups obtained are analyzed using \GAP~\cite{GAP}.
To facilitate the usage of \GAP, we consider the
projective groups $\piproj\CC$ rather than the affine ones
$\piaff\CC$: these groups have finite abelianizations,
hence \GAP\ can easily find their commutants.

In each case, a presentation for $\piproj\CC$ is given by
Corollary~\ref{groups}, with the braid monodromy computed
using the approach of~\cite{degt.kplets} and the skeletons shown in
Figure~\ref{fig.curves}. The reference fiber~$F_0$ is
the trivalent vertex shown in the
figures by a grey dot. Due to Lemma~\ref{infty}, in the presence
of the relation at infinity $\Gr^d=1$, one of the braid relations
$\bm(\Gg_i)=\id$ can be ignored. We omit the relation
resulting from the outermost region of the skeleton (whenever
present).

Each braid relation (resulting from an $m$-gonal region of the
skeleton) has the form $\mu\1\Gs_1^m\mu=\id$ for some
$\mu\in\BG3$.
Given $\Ga,\Gb\in\AA$ and a positive integer~$m$,
introduce the notation
$$
\{\Ga,\Gb\}_m=\cases
(\Ga\Gb)^k(\Gb\Ga)^{-k},&\text{if $m=2k$ is even},\\
\bigl((\Ga\Gb)^k\Ga\bigr)\bigl((\Gb\Ga)^k\Gb\bigr)\1,&
 \text{if $m=2k+1$ is odd}.
\endcases
$$
Then the relation $\mu\1\Gs_1^m\mu=\id$ is equivalent to
$\mu\{\Ga_1,\Ga_2\}_m=1$ or, alternatively, to
$\{\mu(\Ga_1),\mu(\Ga_2)\}_m=1$.
As a special case, the relation $\Gs_2^m=\id$ is equivalent to
$\{\Ga_2,\Ga_3\}_m=1$.

\subsection{The group $\MG(3)$}\label{s.3+3}
The skeleton of~$\CC$ is the tetrahedron shown in
Figure~\ref{fig.curves}(b); the curve is in~$\Sigma_2$ and has
four type~$\tA_2$ singular fibers.
The group $\piproj\CC$ is well known;
various presentations were obtained in a number of papers.
Our approach gives the
relations
$$
\{\Ga_1,\Ga_2\}_3=\{\Ga_2,\Ga_3\}_3=\{\Ga_3,\Ga_2\1\Ga_1\Ga_2\}_3=1,
\quad\Gr^2=1.
$$
One can easily deduce that
$\AM_\CC=(\Lambda\oplus\Lambda)/(t^2-t+1)$; hence
$\Delta_\CC=(t^2-t+1)^2$. In view of Corollary~\ref{larger.group},
one arrives at the following statement.

\corollary\label{cor.Phi6}
If a trigonal curve~$\CC$ admits a uniform
$\GDG(\CG3\oplus\CG3)$-covering, then
$(t^2-t+1)^2\divides|\Delta_\CC(t)$.
\qed
\endcorollary

The converse of Corollary~\ref{cor.Phi6} is discussed in
Remark~\ref{rem.Delta2} below.

Recall that the curve~$\CC$ can be given by the
affine equation
$$
4y^3-(24x^3+3)y+(8x^6+20x^3-1)=0.
$$
The group $\AG4\cong\PSL(2,\FF3)$ of symmetries of the pair
$(\Sigma_2,\CC)$ produces four non-equivalent, although
isomorphic,
torus structures
on~$\CC$, one for each point of
the projective line $\PP^1(\FF3)$ or, equivalently,
for each epimorphism $\GDG(\CG3\oplus\CG3)\onto\DG6$.
In view of Lemma~\ref{torus.invariant}, any other curve admitting
a uniform $\GDG(\CG3\oplus\CG3)$-covering also has at least four
non-equivalent torus structures.

\subsection{The groups $\MG_1(5)$ and $\MG_1(7)$}
The skeletons of the universal curves are shown in
Figures~\ref{fig.curves}(c) and~(d),
respectively. The sets of singular fibers are
$2\bA_4\oplus2\bA_0^*$ in~$\Sigma_2$ for $\MG_1(5)$
and
$3\bA_6\oplus3\bA_0^*$ in~$\Sigma_4$ for $\MG_1(7)$.

These curves do not lead to any surprises: their affine
fundamental
groups are minimal possible, \ie, the semidirect products
$\CG5\rtimes\Z$ and $\CG7\rtimes\Z$, respectively,
the abelianization~$\Z$ acting on
the commutant \via~$-\id$.

\subsection{The group $\MG_1(9)$}
The skeleton of~$\CC$ is shown in
Figure~\ref{fig.curves}(e); the curve is in $\Sigma_6$ and its
singular fibers are $3\tA_8\oplus2\tA_2\oplus3\tA_0^*$. The
relations for $\piproj\CC$ are
$$
\gathered
\Gs_2^i\{\Ga_1,\Ga_2\}_9=
\Gs_2^i\{\Ga_3,(\Ga_2\Ga_1\Ga_2)\1\Ga_1(\Ga_2\Ga_1\Ga_2)\}_1=1,
\quad i=0,1,2,\\
\{\Ga_2,\Ga_3\}_3=1,\qquad
\Gr^6=1.
\endgathered
$$
Denote this group by~$G$. Using the \GAP~\cite{GAP} commands
\vskip\abovedisplayskip
\halign{\qquad\tt#\hss\cr
h := Subgroup(g, [g.1/g.3, g.2/g.3]); \ \# \ N(h)=[g,g]\cr
P := PresentationNormalClosure(g, h);\cr
SimplifyPresentation(P);\cr}
\vskip\belowdisplayskip\noindent
one finds that the commutant $[G,G]$ is generated by three
elements $x_2$, $x_4$, $x_6$ that are subject to the relations
(copying \latin{verbatim})
$$
x_4\1x_2\1x_4x_2=
x_4^{-3}x_2^3=
x_4x_6\1x_4\1x_2x_6x_2\1=
x_2x_4\1x_6x_4x_2\1x_6\1=1,
$$
which simplify to
$$
(x_2x_4\1)^3=[x_2,x_4]=[x_6,x_2x_4\1]=1.
$$
Hence, one has $[G,G]=\<x_2,x_6\>\times\CG3$, with the $\CG3$
factor generated by $x_2x_4\1$.
It follows that $\AM_\CC\cong\Z\oplus\Z\oplus\CG3$ as a group and
$\Delta_\CC(t)=t^2-t+1$.

A similar computation for the
group $\bar G:=G/\Gr^2$ shows that $[\bar G,\bar G]\cong[G,G]$;
since the groups are Hopfian,
the epimorphism $[G,G]\onto[\bar G,\bar G]$ is an isomorphism.
Due to Corollary~\ref{larger.group}, the affine fundamental group
of any trigonal curve admitting a uniform $\DG{18}$-covering factors
to $\piaff\CC\cong(F_2\times\CG3)\rtimes\Z$, where $F_2$ is the
free group on two generators.

\subsection{The group $\MG(5)$}\label{s.Gamma5}
The skeleton is the dodecahedron shown in
Figure~\ref{fig.curves}(f). The curve is in $\Sigma_{10}$ and
has twelve type~$\tA_4$ singular fibers. The relations for
$\piproj\CC$ are
$$
\gathered
\Gs_2^i\{\Ga_1,\Ga_2\}_5=
\Gs_2^i\{\Ga_3,(\Ga_1\Ga_2)\1\Ga_2(\Ga_1\Ga_2)\}_5=1,
\quad i=0,\ldots,4,\\
\{\Ga_2,\Ga_3\}_5=1,\qquad
\Gr^{10}=1.
\endgathered
$$
Denote this group by~$G$, and let
$G_{\AM}=G/\Ga_1^{10}$ and $\bar G=G/\Gr^2$.
Using \GAP~\cite{GAP}, one can easily obtain the following
statements (where $'$ stands for the commutant):
\roster
\item\local{G.1}
$G'\!/G''\cong G_{\AM}'/G_{\AM}''\cong\bigoplus_8\Z$;
\item\local{G.2}
$\bar G'\!/\bar G''\cong\CG5\oplus\CG5$ and
$\bar G''\!/\bar G'''\cong\bigoplus_6\CG5$.
\endroster
Since the Alexander module of~$G_{\AM}$ is obviously a quotient of
both $(\Lambda\oplus\Lambda)/\Phi_{10}$ and $\AM_\CC$,
it follows from~\loccit{G.1}
that $\AM_\CC=(\Lambda\oplus\Lambda)/\Phi_{10}$ and
$\Delta_\CC(t)=\Phi_{10}^2(t)$.
Statement~\loccit{G.2} sheds certain light on the structure
of~$\bar G$; the group appears infinite, but I do not know a
proof: my computer runs out of memory while trying to compute
~$\bar G'''$.

It is also immediate that the epimorphism
$G'\onto\bar G'$ is not an
isomorphism. Hence,
the monodromy group $\BM_\CC(\BG3)$ of the simple universal curve
is a proper subgroup of the pull-back of $\tMG(5)$ in
$\BG3$, \cf. Remark~\ref{rem.simple.proper}.

Due to Corollary~\ref{larger.group}, the fundamental group of any
trigonal curve admitting a uniform $\GDG(\CG5\oplus\CG5)$-covering
factors to~$\bar G$;
due to Corollary~\ref{simple.group}, the group of any
\emph{simple} curve with this property factors to~$G$.

\corollary\label{cor.Phi10}
For a simple irreducible trigonal curve~$\CC$, consider the
following properties\rom:
\roster
\item\local{10.1}
$\piproj\CC$ factors to $\GDG(\CG5\oplus\CG5)$\rom;
\item\local{10.2}
$(t^4-t^3+t^2-t+1)^2\divides|\Delta_\CC(t)$\rom;
\item\local{10.3}
$(t^4-t^3+t^2-t+1)\divides|\Delta_\CC(t)$\rom;
\item\local{10.4}
$\piproj\CC$ factors to $\DG{10}$.
\endroster
Then
\loccit{10.1}$\implies$\loccit{10.2}$\implies$\loccit{10.3}$\implies$\loccit{10.4}.
\endcorollary

\proof
The first implication follows from the computation above and
Corollary \ref{simple.group}, the second one is obvious, and the
last one is given by Lemma~\ref{Delta.2p}.
\endproof

\Remark\label{rem.Delta2}
It can be shown that the converse implication
\iref{cor.Phi10}{10.2}$\implies$\ditto{10.1}, as well as the
converse of Corollary~\ref{cor.Phi6}, also hold; proof will appear
elsewhere. In fact, the Alexander polynomial of a
non-isotrivial trigonal curve
cannot be divisible by $\Phi_{m}^2(t)$ unless $m=2$, $4$, $6$, $8$,
or~$10$.
\endRemark

\Remark
Perturbing the type~$\tA_4$ singular fiber in the central region
of Figure~\ref{fig.curves}(f) to five type~$\tA_0^*$ fibers,
one adds to the presentation above
the relation $\Ga_2=\Ga_3$. Using \GAP~\cite{GAP}, it is easy to
see that the Alexander polynomial of the resulting curve is
$\Phi_{10}(t)$. In particular, both implications
\itemref{cor.Phi10}{10.2}$\implies$\ditto{10.3}$\implies$\ditto{10.4}
in Corollary~\ref{cor.Phi10} are strict.
\endRemark

\section{Further applications\label{S.appl}}

In Subsection~\ref{s.topology}, we relate the maximal uniform
dihedral quotient of $\piaff\CC$
to the homology and the torsion of the Mordell--Weil group of the
covering elliptic surface. In Subsection~\ref{s.Z-splitting}, this
relation is used in the study of another geometric property of the
curve, its so called \emph{$Z$-splitting} sections. Finally, in
Subsections~\ref{s.generalized} and~\ref{proof.plane}, we discuss the
extent to which the results of the paper apply to generalized
trigonal curves and prove Theorem~\ref{th.plane}.

\subsection{The topological interpretation of
uniform dihedral quotients}\label{s.topology}
Fix a trigonal curve $\CC\subset\Sigma_d$. Assume that $d=2k$ is
even. Then there exists a unique double covering $X\to\Sigma_d$
ramified at $\CC+\EE$. Denote by $\tX=\tX_\CC$ the minimal resolution of
singularities of~$X$, and let $\tF_i\subset\tX$, $i=0,\ldots,r$,
be the preimages of the fibers $F_i$ introduced in
Subsection~\ref{s.sections}, \ie, a nonsingular fiber~$F_0$ and
all singular (and possibly more nonsingular) fibers
$F_1,\ldots,F_r$. The pull-backs of~$\CC$ and~$\EE$ in~$\tX$ are
identified with~$\CC$ and~$\EE$ themselves. Note that $\tX$ is a
Jacobian elliptic surface (not necessarily relatively minimal)
and~$\EE$ is its section; one has
$\EE^2=-k$ in~$\tX$.

Denote $\tF_*=\bigcup_{i=1}^r\tF_i$ and $X\scirc=\tX\sminus\tF_*$,
and
let $\tGg_i\in H_1(X\scirc\sminus\EE)$ be the class realized by the
meridian about the proper pull-back of~$F_i$, $i=1,\ldots,r$; this
class can also be realized by any lift of the element~$\Gg_i$ of
any geometric basis $\{\Gg_1,\ldots,\Gg_r\}$ as in
Subsections~\ref{s.sections} and~\ref{s.vanKampen}.

Keeping in mind further applications, fix the following notation:
\Dashes
\dash
$S_\CC\subset H_2(\tX)$ is the subgroup spanned by the
classes of the components of the singular fibers~$\tF_i$,
$i=1,\ldots,r$,
and the section~$\EE$;
\dash
$\chS_\CC=(S_\CC\otimes\Q)\cap H_2(\tX)$ is the primitive hull
of~$S_\CC$;
\dash
$\CK_\CC=\chS_\CC/S_\CC$;
\dash
$\CQ_\CC$ is the group $\CQ=\CH/H$
given by Corollary~\ref{cor.CQ}.
\endDashes
Thus, $S_\CC$ is the image of the inclusion homomorphism
$\inj_*\:H_2(\tF_*\cup\EE)\to H_2(\tX)$, and $\CK_\CC$ measures
the imprimitivity of this image. The subgroup
$\chS_\CC\subset H_2(\tX)$ can be regarded as the `minimal'
Neron--Severi group, \ie, the Neron--Severi group of the
surface~$\tX$ obtained from a curve generic in its equisingular
fiberwise deformation family.

\lemma\label{H1X}
There are canonical isomorphisms
$$
H_1(X\scirc)=H_1(X\scirc\sminus\EE)=\CQ_\CC\oplus\CZ,
$$
where $\CZ$ is the free abelian group
$\bigoplus_{i=1}^r\Z\tGg_i/\sum_{i=1}^r\tGg_i$. The deck
translation acts \via\ $\id$ and~$-\id$ on~$\CZ$ and~$\CQ_\CC$,
respectively.
\endlemma

\proof
The statement is immediate both algebraically, using
Theorem~\ref{th.vanKampen}, and topologically, computing the
$1$-homology of
the locally trivial fibration $X\scirc\to\BB\scirc$ or
$X\scirc\sminus\EE\to\BB\scirc$
using its monodromy.
In both approaches,
crucial is the fact that the monodromy in the homology
$$
\CH=H_1(\tF)=H_1(\tF\sminus\EE)=
 \Ker{\deg_2}/\<\Ga_1^2,\Ga_2^2,\Ga_3^2\>,
 \eqtag\label{eq.CH}
$$
see Subsection~\ref{s.quotients}, of the pull-back
of a
reference fiber~$F$ (the homological invariant of
the elliptic surface~$\tX$)
is the $\tMG$-valued reduction of the braid
monodromy of~$\CC$, see Lemma~\ref{induced.action}.
Note also that $\Gr^2$ projects to
$0\in H_1(\tF\sminus\EE)$, hence the relation at infinity becomes
redundant.
\endproof

\corollary\label{Q=K}
Unless $\CC$ is trivial, there is an isomorphism
$\CQ_\CC=\Ext(\CK_\CC,\Z)$.
\endcorollary

\proof
Since $\CC$ is nontrivial, so is~$\tX$ and $H_1(\tX)=0$.
From the
Poincar\'e--Lefschetz duality and the exact sequence of pair
$(\tX,\tF_*)$ one has
$$
H_1(\tX\scirc\sminus\EE)=
 \Coker[\inj^*\:H^2(\tX)\to H^2(\tF_*\cup\EE)].
$$
Since all homology groups involved are torsion free, the
torsion of the last
cokernel equals $\Ext(\CK_\CC,\Z)$.
On the other hand, due to Lemma~\ref{H1X} and the fact that
$\CQ_\CC$ is a torsion group, see Theorem~\ref{th.D},
one has
$\CQ_\CC=\Tors H_1(\tX\scirc\sminus\EE)$.
\endproof

\Remark\label{rem.Q=K}
Since both~$\CQ_\CC$ and $\CK_\CC$ are finite groups, the
isomorphism given by Corollary~\ref{Q=K} can be rewritten in the
form $\CQ_\CC=\Hom(\CK_\CC,\Q/\Z)$ or, in other words, there is a
nonsingular pairing $\CQ_\CC\otimes\CK_\CC\to\Q/\Z$. Hence,
conversely, one has
$\CK_\CC=\Hom(\CQ_\CC,\Q/Z)$.
\endRemark

\corollary\label{Tors.MW}
Assume that the trigonal curve $\CC\subset\Sigma_{2k}$ is nontrivial.
Then the map sending a section of~$\tX$ to its homology class establishes
isomorphisms $\Tors\MW(\tX)=\CK_\CC=\Hom(\CQ_\CC,\Q/Z)$,
where $\MW$ is the Mordell--Weil group.
\endcorollary

\proof
According to T.~Shioda~\cite{Shioda}, since $H_1(\tX)=0$,
the map
above
establishes an isomorphism
$\MW(\tX)=\NS(\tX)\!/S_\CC$ and, since $\NS(\tX)$ is
primitive in $H_2(\tX)$,
one has
$\Tors\MW(\tX)=\Tors(H_2(\tX)\!/S_\CC)=\CK_\CC$.
It remains to apply Corollary~\ref{Q=K}.
\endproof

As a further consequence, we obtain the following known
theorem.

\theorem[Theorem \rm(D.~A.~Cox, W.~R.~Parry~\cite{Cox})]\label{th.MW}
The torsion of the Mordell--Weil
group of a non-trivial
Jacobian elliptic surface over a rational base is
isomorphic to one of the groups~$\CQ$ given by Theorem~\ref{th.D}.
\qed
\endtheorem

\subsection{$Z$-splitting sections}\label{s.Z-splitting}
The notion of $Z$-splitting curves was introduced by
I.~Shimada in~\cite{Shimada1}; we remind the definition
of~\cite{Shimada1}, modifying it
for the case of trigonal curves and splitting sections.
Splitting, although not $Z$-splitting, sections for trigonal curves
were also studied in~\cite{Shioda2} and~\cite{Tokunaga.splitting}.

We keep the setting and the notation of
Subsection~\ref{s.topology} and assume that the trigonal
curve $\CC\subset\Sigma_d$, $d=2k$, is nontrivial.
Now, the word `section' stands for a
\emph{holomorphic} section of the ruling $\Sigma_d\to\Cp1$.

\definition\label{def.splitting}
A section~$L\subset\Sigma_d$ is called \emph{splitting} for~$\CC$
if the proper transform of~$L$ in~$\tX$ splits into two
distinct components $L_+$ and~$L_-$.
A splitting section is called \emph{pre-$Z$-splitting} if
$[L_\pm]\in\chS_\CC$; in this case, the order of $[L_\pm]$ in the
finite group $\CK_\CC=\chS_\CC/S_\CC$
is called the \emph{class order} of~$L$.
Finally, a pre-$Z$-splitting section is called
\emph{$Z$-splitting} if $[L_+]\ne[L_-]$.
\enddefinition

Informally, $Z$-splitting sections are those that remain stable
under equisingular fiberwise deformations of the curve, as one
always has $\chS_\CC\subset\NS(\tX)$. Note though that, in the
exceptional case $\CC\subset\Sigma_2$, when \smash{$\tX$} is rational, the
curve may have other stable splitting sections, which are not
$Z$-splitting.

Clearly, the components $L_\pm$ of the pull-back of a splitting
section~$L$ form a pair of sections of
the Jacobian elliptic surface~$\tX$;
they are interchanged by the deck
translation. Conversely, any section of~$\tX$ that is \emph{not}
deck translation invariant projects to a splitting section
for~$\CC$. Since
a section of~$\tX$ is uniquely determined by its homology class,
see~\cite{Shioda},
the last condition $[L_+]\ne[L_-]$ in
Definition~\ref{def.splitting} is redundant and any
pre-$Z$-splitting section is $Z$-splitting. Finally, according to
Corollary~\ref{Tors.MW}, a splitting section is $Z$-splitting
if and only if
the components~$L_\pm$ are torsion elements of $\MW(\tX)$,
and the class order of~$L$ equals the order of
$L_\pm$ in $\MW(\tX)$.
Combining, one arrives at the following statement.

\theorem\label{th.splitting}
There is a canonical one-to-one correspondence between the set of
$Z$-splitting sections of a nontrivial trigonal curve~$\CC$ and
the set of unordered pairs of distinct opposite elements of
$\CK_\CC=\Hom(\CQ_\CC,\Q/\Z)$.
\qed
\endtheorem

It follows that the number of $Z$-splitting sections
for~$\CC$
is $\frac12(\ls|\CQ_\CC|-\ls|\CQ_\CC\otimes\CG2|)$.
In particular, an irreducible curve may have
zero to four or twelve
such sections.
If $\CC$ is reducible, the group $\Hom(\CQ_\CC;\Q/\Z)$ has one or
three elements of order~$2$, see Theorem~\ref{th.reducible}; they
represent deck translation invariant sections of~$\tX$.
Clearly, these
sections are the components of~$\CC$ that are sections
of~$\Sigma_d$.

Corollary~\ref{universal.simple} gives one a very explicit
description of $Z$-splitting sections: each such section~$L$
is a
pull-back and/or Nagata transform of a $Z$-splitting section~$L'$
for an
appropriate universal curve; we will say that $L$ is
\emph{induced} from~$L'$.
Indeed, all sections obtained in
this way are $Z$-splitting, and their number equals that given by
Theorem~\ref{th.splitting}.
Furthermore,
an order~$m$ element of $\CK_\CC=\Hom(\CQ_\CC,\Q/Z)$ can be
regarded as
an epimorphism $\CQ_\CC\to\CG{m}$;
hence
the corresponding $Z$-splitting section of class order~$m$
is induced from a section for the universal curve~$C_m$
corresponding to~$\MG_1(m)$.
Thus, there is a finite list of
essentially distinct $Z$-splitting sections.

Any $Z$-splitting
section of class order~$3$ is the section $y+a_2(x)=0$
appearing in the corresponding torus structure,
see~\eqref{eq.torus}. For proof,
it suffices to check that
the section $y=0$
is $Z$-splitting
for the universal curve
$C_3=\{y^3+(y+x)^2=0\}$
corresponding to $\MG_1(3)$,
see~\ref{ss.3} and Figure~\ref{fig.curves}(a).

The curve~$C_5\subset\Sigma_2$, see Figure~\ref{fig.curves}(c),
can be given by
the Weierstra{\ss} equation \eqref{eq.W} with
$$
\gathered
-4p(x)=x^4-12x^3+14x^2+12x+1,\\
8q(x)=(x^2+1)(x^4-18x^3+74x^2+18x+1),
\endgathered
$$
and any $Z$-splitting section of class order~$5$ is induced from
one of the two
sections for~$C_5$; they are given by $2y=x^2\pm6x+1$. Each of
these
sections passes through both type~$\bA_4$ singular points of~$C_5$
and is tangent to~$C_5$ at one of these points.
Note that the two
$Z$-splitting
sections are interchanged by the involutive symmetry
$(x,y)\mapsto(1/x,y/x^2)$ of~$C_5$, constituting a single
isomorphism class.

The
curve~$C_4\subset\Sigma_2$, see
Figure~\ref{fig.curves.r}(c),
is given by $(y^2-4x)(y-x-1)=0$, and its
$Z$-splitting section of class order~$4$ is $y=2$; it passes
through the type~$\bA_3$ point of~$C_4$ (over $x=1$)
and is tangent to the cusp
at its type~$\bD_5$ point (over $x=\infty$).

The curve~$C_6\subset\Sigma_2$, see
Figure~\ref{fig.curves.r}(i), is given by $(y^2-x)(y-l(x))=0$,
where
$16l=-16x^2+24x+3$.
Its $Z$-splitting section of class order~$3$ (the one arising from
the torus structure) is $4y=4x+1$; it passes through the
type~$\bA_2$ singular point of~$C_6$ (over $x=\infty$)
and is tangent to~$C_6$ at
its type~$\bA_5$ singular point (over $x=1/4$).
The $Z$-splitting section
of class order~$6$ is $4y=-4x+3$;
it
passes through all
three singular points of~$C_6$, including the type~$\bA_1$ point
over $x=9/4$.

For the remaining class orders~$7$, $8$, $9$, and~$10$, I do not
know explicit equations of the corresponding curves. Using
Nikulin's theory of discriminant forms, it is not difficult do
determine how a $Z$-splitting section should pass through the
singularities of~$C_m$; I leave this computation to the reader.

\subsection{Generalized trigonal curves}\label{s.generalized}
By a \emph{generalized trigonal curve} we mean a reduced curve
$\CC\subset\Sigma_d$ with the following properties:
\roster
\item
$\CC$ does not contain a fiber of~$\Sigma_d$ as a component;
\item
$\CC$ intersects each fiber of~$\Sigma_d$ at three points.
\endroster
Thus, unlike genuine trigonal curves defined in
Subsection~\ref{s.curves}, a generalized trigonal curve is allowed
to meet the exceptional section~$\EE$. A \emph{singular fiber} of
a generalized trigonal curve $\CC\subset\Sigma_d$ is a fiber
of~$\Sigma_d$ intersecting $\CC+\EE$ geometrically at fewer than
four points.

As in the case of genuine trigonal curves, see
Subsection~\ref{s.results},
we will consider
both the
projective fundamental group
$\piproj\CC:=\pi_1(\Sigma_d\sminus(\CC\cup\EE))$
and
the affine group
$\piaff\CC:=\pi_1(\Sigma_d\sminus(\CC\cup\EE\cup F))$, where
$F\subset\Sigma_d$
is a generic nonsingular fiber.

By a sequence of positive Nagata transformations, any generalized
trigonal curve $\CC\subset\Sigma:=\Sigma_d$ can be transformed to a
genuine trigonal curve $\CC'\subset\Sigma':=\Sigma_{d'}$
for some $d'\ge d$. The latter is called a \emph{trigonal model}
of~$\CC$; it is defined up to Nagata equivalence of trigonal
curves. Let $F_i$, $i=1,\ldots,r$, be all fibers that are singular
to either~$\CC$ or~$\CC'$, and let~$F_0$ be a common nonsingular
fiber.
(Recall that we identify fibers of the
ruling and their projections to the base. Hence, we can also
identify fibers of~$\Sigma$ and those of~$\Sigma'$.) The Nagata
transformations used establish a biholomorphism
$$
\Sigma\sminus(\CC\cup\EE\cup\bigcup_{i=0}^rF_r))\cong
 \Sigma'\sminus(\CC'\cup\EE\cup\bigcup_{i=0}^rF_r)),
$$
hence an isomorphism of the respective fundamental groups. Thus,
the presentation given by Theorem~\ref{th.vanKampen} still holds
for a generalized trigonal curve~$\CC$, assuming that the braid
monodromy~$\bm$ is defined using a trigonal model~$\CC'$.

\Remark
The braid monodromy~$\bm$ does depend on the trigonal
model~$\CC'$. Nevertheless, the presentations given by
Theorem~\ref{th.vanKampen} result in isomorphic groups: they
differ by the lifts~$\tGg_i$. A positive Nagata transformation
contracting fiber~$F_i$ multiplies~$\bm_i$ by~$(\Gs_2\Gs_1)^3$, and
the new lift of~$\Gg_i$ is $\tGg_i\Gr\1$.
\endRemark

Fix a trigonal model~$\CC'$ and a geometric basis
$\{\Gg_1,\ldots,\Gg_r\}$ as in Subsection~\ref{s.sections}.
Then, to each basis element~$\Gg_i$ one can assign a \emph{slope}
$\slope_i\in\AA$, see~\cite{degt.dessin}. (Note that the slope
depends on both the curve~$\CC'$ and generator~$\Gg_i$.)
According to~\cite{degt.dessin}, patching back in a fiber~$F_i$
results in an extra relation $\tGg_i\slope_i=1$. Finally,
Corollary~\ref{groups} takes the following form.

\corollary\label{groups.generalized}
For a generalized trigonal curve $\CC$,
geometric presentations of the affine and projective fundamental
groups are as follows\rom:
$$
\align
\piaff\CC&=\AA/
 \<\text{$\slope_i\Ga\slope_i\1=\bm_i(\Ga)$, \ $i=1,\ldots,r$, \ $\Ga\in\AA$}\>,\\
\piproj\CC&=\piaff\CC\!/\<\slope_r\ldots\slope_1\Gr^{d'}\>,
\endalign
$$
where the monodromies $\bm_i\in\BG3$ and slopes $\slope_i\in\AA$,
$i=1,\ldots,r$,
are defined using a trigonal model $\CC'\subset\Sigma_{d'}$.
\qed
\endcorollary

It is shown in~\cite{degt.dessin} that $\piaff\CC$ is an extension
of~$\piproj\CC$ by a central infinite cyclic subgroup. In
particular, the commutants of the two groups are isomorphic.

\definition\label{def.even}
A sequence of Nagata transformations converting a generalized
trigonal curve~$\CC$ to a trigonal model~$\CC'$ is called
\emph{even} if all slopes $\slope_i$ are of even degree. (This
property is independent of the choice of a basis
$\{\Gg_1,\ldots,\Gg_r\}$, see~\cite{degt.dessin}.)
If this is the case, $\CC'$ is called
an \emph{even trigonal model} of~$\CC$.
\enddefinition

It is easy to see that a trigonal model of~$\CC$ is even if and
only if it differs from the divisorial transform of~$\CC$ by an
even divisor. (For example, one can use the computation of local
slopes found in~\cite{degt.dessin}.) A trigonal model of a genuine
trigonal curve~$\CC$ is even if and only if it is $2$-Nagata
equivalent to~$\CC$; in particular,
all even trigonal models of any curve~$\CC$
are $2$-Nagata equivalent.
Any generalized trigonal curve admits an even
trigonal model: each slope $\slope_i$ of odd degree can be
corrected by a single Nagata transformation in the corresponding
fiber: it multiplies~$\slope_i$ by~$\Gr$.

Consider the maximal uniform dihedral quotient $\AA\to\GDG(\CH)$,
see Subsection~\ref{s.quotients}, and, for an element $\Ga\in\AA$ of
even degree, denote by $\tGa\in\CH$ its projection to~$\CH$.

Fix a generalized trigonal curve~$\CC$ and pick an
even trigonal model~$\CC'$ of~$\CC$. Denote by~$\slope_i$,
$i=1,\ldots,r$, the slopes and let $\piaff{\CC'}\onto\GDG(\CQ')$,
$\CQ'=\CH/H$, be the maximal uniform dihedral quotient, see
Corollary~\ref{cor.CQ}.

\theorem\label{th.generalized}
In the notation above, the maximal uniform dihedral quotient
of the affine group $\piaff\CC$ is $\GDG(\CQ\saff)$, where
$\CQ\saff=\CQ'\!/\<2\tslope_i\>$, $i=1,\ldots,r$.
If, in addition, $\CC'\subset\Sigma_{d'}$ with $d'=2k$ even, then
the maximal uniform dihedral quotient of
the projective group $\piproj\CC$ is
$\GDG(\CQ\sproj)$, where
$\CQ\sproj=\CQ\saff\!/\<\tslope_1+\ldots+\tslope_r\>$.
\endtheorem

\proof
The statement follows from Corollary~\ref{groups.generalized}. For
the projective group $\piproj\CC$, it suffices to consider the
braid relations $\slope_i\Ga\slope_i\1=\bm_i(\Ga)$ for $\Ga\in\AA$
an element of even degree and for $\Ga=\Gr$. In the former case,
one has $\tslope_i\tGa\tslope_i\1=\tilde\bm_i(\tGa)$, where $\tilde\bm_i$
is the induced automorphism of~$\CH$. Since $\CH$ is
abelian, this implies the relations $\tilde\bm_i=\id$, which define~$\CQ'$,
see Corollary~\ref{cor.CQ}. In the latter case, since
$\bm_i(\Gr)=\Gr$ and $\Gr\1\tslope_i\Gr=-\tslope_i$ in~$\CH$ (the
definition of dihedral groups; recall that $\deg_2\Gr=1$), one has
$2\tslope_i=0$. For the projective group, since $d'$ is even and
$(\Gr^2)\sptilde=0$, the relation at infinity results in
$\tslope_1+\ldots+\tslope_r=0$.
\endproof

\Remark
If $d'$ is odd in Theorem~\ref{th.generalized}, then the maximal
uniform cyclic quotient of $\piproj\CC$ is of odd order and
$\piproj\CC$ has no dihedral quotients.
\endRemark

A generalized trigonal curve is called \emph{trivial} if its even
trigonal models are trivial.

\corollary\label{cor.generalized}
Theorems~\ref{th.D} and~\ref{th.reducible} and
Corollary~\ref{cor.irreducible} hold literally
for generalized trigonal
curves.
\qed
\endcorollary

\subsection{Proof of Theorem~\ref{th.plane}}\label{proof.plane}
Blow up a singular point of~$\DD$ of multiplicity $(\deg\DD-3)$.
The result is the Hirzebruch surface~$\Sigma_1$, and the proper
transform of~$\DD$ is a generalized trigonal curve
$\CC\subset\Sigma_1$. The projection $\Sigma_1\to\Cp2$ restricts
to a biholomorphism
$\Sigma_1\sminus(\CC\cup\EE)\cong\Cp2\sminus\DD$, hence inducing
an isomorphism
$\piproj\CC=\pi_1(\Cp2\sminus\DD)$, and the statement
follows from Corollary~\ref{cor.generalized}.
\qed


\example\label{ex.Klein}
Take for $\DD\subset\Cp2$ the three cuspidal quartic. As is well
known, $\pi_1(\Cp2\sminus\DD)$ is the alternating group~$\AG4$,
and the maximal dihedral quotient of this group is $\DG6$. Any
smooth point $P\in\DD$ can be taken for
the distinguished `singular' point of
multiplicity $(\deg\DD-3)=1$. Hence the curve satisfies the
hypotheses of Theorem~\ref{th.plane}. However, in spite of the
existence of a $\DG6$-covering, Statement~\iref{th.Oka}{Oka.D},
none of
the other statements of Theorem~\ref{th.Oka}
holds for~$\DD$.

If the point $P$ to be blown up is generic, the minimal trigonal
model of~$\DD$ is a curve $\CC'\subset\Sigma_2$ with the set of
singular fibers $3\tA_2\oplus\tA_1\oplus\tA_0^*$. Up to $2$-Nagata
equivalence, $\CC'$ is induced from the universal curve
corresponding to $\tMG_1(3)$, see~\ref{ss.3}.
If $P$ is in the double tangent to~$\DD$, the trigonal model is
the universal curve corresponding to $\tMG(3)$, see
Subsection~\ref{s.3+3}.
\endexample


\subsection{Speculation~\ref{spec} for generalized trigonal
curves}\label{s.spec}
We conclude with a few remarks concerning extending
Statements~\iref{spec}{spec.universal}--\ditto{spec.larger} to
generalized trigonal curves and, in particular, to plane curves
with a singular point of multiplicity $({\deg}-3)$ (which can be
regarded as generalized trigonal curves in~$\Sigma_1$). We keep
the notation of Theorem~\ref{th.generalized}: $\CC\subset\Sigma_d$
is a generalized trigonal curve, $\CC'\subset\Sigma_{d'}$ is an even
trigonal model of~$\CC$, and $\CQ'$ and
$\CQ\saff=\CQ'\!/\<2\tslope_i\>$,
$i=1,\ldots,r$, are the
abelian groups appearing in the maximal uniform
dihedral quotients
$\piaff{\CC'}\onto\GDG(\CQ')$ and
$\piaff\CC\onto\GDG(\CQ\saff)$.
One has an epimorphism $\CQ'\onto\CQ\saff$; in general, it may be
proper.

\subsubsection{Relation to universal curves}
The existence of a uniform dihedral quotient
$\piaff\CC\onto\GDG(\CG{m}\oplus\CG{n})$, $m\divides|n$,
still implies that, up to even Nagata equivalence, $\CC$ is
induced from the universal curve corresponding to $\tMG_m(n)$.
However, the converse is no longer true and
Corollaries~\ref{larger.group} and~\ref{simple.group} may fail,
\cf.~Example~\ref{ex.Klein}.
The Alexander polynomial may also change.

\subsubsection{Torus structures}
Torus structures do not survive in~$\CC$: in general, one can only
assert that the sum of~$\CC$ with a certain linear combination of
fibers is a curve of torus type (in the obvious sense), see
Remark~\ref{rem.torus}. For example, the three cuspidal quartic~$\DD$
in Example~\ref{ex.Klein} is not of torus type but, for any
tangent~$L$ to~$\DD$, the sextic $\DD+2L$ is of torus type.
(After the tangency point is blown up, $L$ becomes the fiber
contracted by the Nagata transformation.)
If $L$
is the (only) double tangent to~$\DD$,
then $\DD+2L$ admits four non-equivalent
torus structures, as in this case the trigonal model of~$\DD$ is
the universal curve corresponding to~$\tMG(3)$, see
Subsection~\ref{s.3+3}.

\subsubsection{$Z$-splitting sections}
Clearly, the $Z$-splitting sections for~$\CC$ are precisely the
transforms of those for~$\CC'$, the correspondence preserving the
class order; in fact, under an even sequence of Nagata
transformations, the corresponding double coverings $\tX$
and~$\tX'$ differ at most by a few blow-ups at nonsingular points.
Hence, the $Z$-splitting sections of~$\CC$ are enumerated by pairs
of distinct opposite elements of $\Hom(\CQ',\Q/\Z)$, see
Theorem~\ref{th.splitting},
but, in general, not of $\Hom(\CQ\saff,\Q/\Z)$.
It may happen that $\CC$ has $Z$-splitting sections but
$\CQ\saff=0$, \ie, $\piaff\CC$ does not admit
a single dihedral quotient.

\refstyle{C}

\enddocument